\documentclass[11pt]{article}
\usepackage{amscd,amsmath,amsxtra,amssymb,theorem,latexsym,amsfonts,color,graphics}
\usepackage{ulem}

\newtheorem{Thm}{Theorem}[section]
\newtheorem{Lemma}[Thm]{Lemma}
\newtheorem{Prop}[Thm]{Proposition}

\newtheorem{Notation}[Thm]{Notation}

\newcommand{\R}{{\mathbb R}}

\newcommand{\C}{{\mathbb C}}
\newcommand{\Z}{{\mathbb Z}}
\newcommand{\N}{{\mathbb N}}

\newcommand{\id}{{\rm id}}

\newcommand{\bfA}{{\mathbf A}}
\newcommand{\bfB}{{\mathbf B}}
\newcommand{\bfF}{{\mathbf F}}
\newcommand{\bfH}{{\mathbf H}}

\newcommand{\bfg}{{\mathbf g}}

\newcommand{\bfw}{{\mathbf w}}
\newcommand{\bfx}{{\mathbf x}}
\newcommand{\bfy}{{\mathbf y}}

\newcommand{\calE}{{\cal E}}
\newcommand{\calF}{{\cal F}}
\newcommand{\calG}{{\cal G}}
\newcommand{\calH}{{\cal H}}
\newcommand{\calI}{{\cal I}}
\newcommand{\calK}{{\cal K}}
\newcommand{\calL}{{\cal L}}

\newcommand{\calO}{{\cal O}}

\newcommand{\calT}{{\cal T}}
\newcommand{\calU}{{\cal U}}

\newcommand{\calX}{{\cal X}}
\newcommand{\calY}{{\cal Y}}
\newcommand{\calZ}{{\cal Z}}

\newcommand{\tOmega}{{\tilde \Omega}}

\newcommand{\sfA}{{\mathsf{A}}}
\newcommand{\sfV}{{\mathsf{V}}}

\newcommand{\rmh}{{\rm h}}

\begin{document}

\title{An example of Newton's method for an equation in Gevrey series}
\author{Alexander GETMANENKO \\
\footnotesize Kavli IPMU, University of Tokyo, Japan; \\ \footnotesize Mathematics Institute of Jussieu, UPMC, Paris, France; \\ \footnotesize  {\tt getmanenko@math.jussieu.fr}.}
\maketitle

\begin{abstract} In the context of complex WKB analysis, we discuss a one-dimensional Schr\"odinger equation
$$ -h^2\partial_x^2 f(x,h) + [Q(x)+hQ_1(x,h)]f(x,h) =0, \ \ \ h\to 0, $$
where $Q(x)$, $Q_1(x,h)$ are analytic near the origin $x=0$, $Q(0)=0$, and $Q_1(x,h)$ is a factorially  divergent power series in $h$. We show that there is a change of independent variable $y=y(x,h)$, analytic near $x=0$ and factorially divergent with respect to $h$, that transforms the above Schr\"odinger equation to a canonical form. The proof goes by reduction to a mildly nonlinear equation on $y(x,h)$ and by solving it using an appropriately modified Newton's method of tangents.
\end{abstract}


\section{Introduction}

{\bf \large Context.} 
The stationary Schr\"odinger equations for an anharmonic oscillator
\begin{equation}  -h^2 \partial_x^2 f(x,h) + Q(x) f(x,h) = 0, \label{BeforeDef} \end{equation}
where $h\to 0+$, $f(x,h)$ is an unknown function, and $V(x)$ is, say, a polynomial with $$ Q(0)=Q'(0)=0\ne Q''(0), $$  is one of the basic problems in Mathematical Physics. It is easy to find formal WKB solutions of \eqref{BeforeDef} by substituting an ansatz \begin{equation} f(x,h) \sim e^{S(x)/h}  (a_0(x)+a_1(x)h+a_2(x)h^2+...) \label{WKBansatz} \end{equation} and recursively solving for new unknown functions $S(x),a_0(x), a_1(x),.. $, yet analytic properties of this expansion are quite subtle. 
Even in simplest examples, the terms of the series $a_0(x)+a_1(x)h+a_2(x)h^2+...$ have singularity at the origin and diverge factorially (i.e., are {\it Gevrey series}, sec. \ref{GevreyNorms}) away from zeros of $Q(x)$; a lot of effort is being put into proving Borel summability divergent series of this origin.

The singular behavior at the origin of the coefficients $a_j(x)$ should not be all that uncontrollable: the closer $x$ gets to $0$, the more the equation \eqref{BeforeDef} resembles the Schr\"odinger equation of a harmonic oscillator whose eigenfunctions can be expressed in terms of well-known cylindro-parabolic functions. 

And in fact, in ~\cite{AKT} it is shown that one can find a change of coordinates $y=y(x,h)$, where $y$ is a no more than factorially divergent power series in $h$ with coefficients holomorphic functions of $x$ {\it in a full neighborhood of the origin}  that reduces \eqref{BeforeDef} to a harmonic oscillator. The proof goes by reducing the problem to a mildly nonlinear equation (case $M=2$ and $Q_1(x,h)=0$ of our \eqref{ma17f2}, see also \eqref{MildlyNL}),   recursively solving it for $E_{0,n},T_n(z)$ where $E_0=\sum_{n\ge 0} E_{0,n} h^n$, $T=\sum_{n\ge 0} T_n(z) h^n$ and analyzing the growth of $|E_{0,n}|$, $|T_n(z)|$. Much more is contained in the rich and beautiful paper ~\cite{AKT}, but here we will concentrate just on this aspect.

{\bf \large Motivation.} 
 Our interest in the Witten Laplacian ~\cite{G} led us to study the paper ~\cite{AKT} with the task to extend its results to equations of slightly more general type than \eqref{BeforeDef}, namely:
\begin{equation}  -h^2 \partial_x^2 f(x,h) + (Q(x)+hQ_1(x,h)) f(x,h) = 0, \label{AfterDef} \end{equation}
where $Q_1(x,h)=\sum_{n\ge 0} Q_{1,n}(x) h^n$ is a factorially divergent series in $h$, $Q(0)=Q'(0)=0\ne Q''(0)$; in particular, to construct a coordinate change $y=y(x,h)$ in a full neighborhood of $x=0$ that would deform \eqref{AfterDef} to a harmonic oscillator. By the time we asked all {\it why} and {\it why not} questions about the original proof of ~\cite{AKT}, we ended up with their proof repackaged in the form of  Newton's method of tangents in the spaces of factorially divergent series. This version gives us the desired generalization automatically. Just as automatically we have obtained a generalization to an arbitrary order of zero of $Q(x)$. 

The {\bf \large main result} is: 

\begin{Thm} \label{MAINTH} Suppose that in the equation \eqref{AfterDef}: \\ 
a) all functions $Q(x)$, $Q_{1,n}(x)$, $n\ge 0$, are analytic on a common neighborhood $U\subset \C$ of $x=0$; \\ b)  $Q(0)=Q'(0)=...=Q^{(M-1)}(0)=0\ne Q^{(M)}(0)$, $M\ge 1$, \\ c)  there are constants $\tau$, $C_0$ such that $\sup_{x\in U} |Q_{1,n}(x)| \le C_0 \tau^n n!$. \\ Then there exist: \\ 1) series $y=y(x,h) = \sum_{n\ge 0} y_n(x) h^n$, where $y_n(x)$ are holomorphic functions on a common neighborhood $U'\subset \C$ of $x=0$, \\ 2)  series $E_j(h) = \sum_{n\ge 0} E_{j,n} h^n$, $E_{j,n}\in \C$, for $0\le j\le M-2$, 
\\ 3) series $\psi=\psi(x,h) = \sum_{n\ge 0} \psi_n(x) h^n$, $\psi_n(x)$ analytic on $U'$; 
\\ 4)  constants $C'_0$, $\tau'$ \\ so that $$ \sup_{x\in U'} |y_n(x)| \le C'_0 {\tau'}^n n! ,  \ \  \ \sup_{x\in U'} |\psi_n(x)| \le C'_0 {\tau'}^n n!, \ \ \  |E_{j,n}|\le C'_0 {\tau'}^n n!, $$
 and the change of independent variable $y=y(x,h)$ followed by rescaling of the unknown function $f(x,h)=\psi(x,h) g(y(x),h)$ in \eqref{AfterDef} transforms this equation to
$$ -h^2\partial_y^2 g(y,h) + (hE_0(h) + hE_1(h)y+...+hE_{M-2}(h) y^{M-2} - \frac{y^M}{4}) g(y,h) = 0 .$$
\end{Thm}

\vskip2.5pc

{\bf \large Plan of the paper.} Section \ref{MthOrderFormal} is devoted to formulating our problem in terms of a mildly nonlinear equation \eqref{ma17f2}. Sections \ref{GevreyNorms} and \ref{FuncAn} introduce definitions of functional spaces used in the rest of the paper. The section \ref{secLA3} discusses inversion of the linear operator which is the dominant member of  equation \eqref{ma17f2}. Section \ref{NMblah} contains no precise definitions or results but explains the intuition based on the Newton's method behind our  argument. Section \ref{Gcalc} contains quite general results about analytic operations with factorially divergent expansions. The main result is proven in section \ref{MainPf}.

The {\bf \large contribution} of this paper is, in our understanding, as follows.  Theorem \ref{MAINTH} is a new generalization of  ~\cite[Theorem A.1]{AKT}. In the vast literature on Gevrey solutions of nonlinear differential equations we have not seen any applications of the Newton's method of tangents, so this way of approaching the problem is probably also new. The intermediate result of section  \ref{Gcalc} should be well known, but at the moment we do not have a specific reference.

\section{Formal reduction to the normal form.} \label{MthOrderFormal}


Following the ideas of ~\cite[{\S}1 and Appendix A]{AKT} whose treatment corresponds to ours in case $M=2$, $Q_1=0$, we consider an equation  
\begin{equation} - h^2 \frac{ \partial^2 f(x,h)}{\partial x^2} + (Q(x)+hQ_1(x,h))f(x,h) = 0,  \label{ma17e1} \end{equation}
where $Q(x)$ is an analytic function in a neighborhood $U\subset \C$ of the origin, $Q(x)=cx^M + O(x^{M+1})$ as $x\to 0$ for some $c\ne 0$, and $Q_1(x,h)=\sum_{n\ge 0} h^n Q_{1,n}(x)$ is a formal (for now) power series in $h$ with coefficients $Q_{1,n}(x)$ analytic functions on a common neighborhood of $x=0$. 

 We seek to simplify \eqref{ma17e1} by an invertible change of independent variable $y=y(x,h)=\sum y_j(x) h^j$, which is (for now) a formal power series in $h$ with coefficients $y_j(x)$ holomorphic functions in a common neighborhood of $0$. Then \eqref{ma17e1} becomes
\begin{equation}  -h^2 \left(  \frac{d^2f}{dy^2} -\frac{d^2x}{dy^2} \left( \frac{dx}{dy}\right)^{-1}  \frac{df}{dy} \right)  +\left(\frac{dx}{dy} \right)^{2}(Q(x)+hQ_1(x,h))f(x) = 0. \label{ma17e2} \end{equation}

Let us further replace the unknown function $f(x,h)$ with a product $\psi(x,h) g(x,h)$, where $g(x,h)$ will be the new unknown function and $\psi$ will be conveniently chosen:
\small
$$  -h^2 \left( -\frac{d^2x}{dy^2} \left( \frac{dx}{dy}\right)^{-1}  
\left[ \psi \frac{dg}{dy} + g\frac{d\psi}{dy} \right] +  \left[ \psi \frac{d^2g}{dy^2} + 2\frac{d\psi}{dy} \frac{dg}{dy} + g \frac{d^2\psi}{dy^2} \right] \right)  +\left(\frac{dx}{dy} \right)^{2}(Q(x)+hQ_1(x,h))\psi(x)g(x) = 0. $$
\normalsize
The following choice eliminates the summands containing $\frac{\partial g}{\partial y}$: 
\begin{equation}  \psi(x,h) = \exp \left\{ \frac{1}{2} \int \frac{d^2 x}{dy^2} \left( \frac{dx}{dy} \right)^{-1} dy \right\}. \label{factorpsi} \end{equation} 
The equation \eqref{ma17e1} is thus reduced to
\begin{equation}  -h^2 \frac{d^2g}{dy^2}    +\left\{ \left(\frac{dx}{dy} \right)^{2}(Q(x)+hQ_1(x,h)) + h^2 \frac{d^2x}{dy^2} \left( \frac{dx}{dy}\right)^{-1}  \frac{(\partial\psi/\partial y)}{\psi} 
- h^2  \frac{(\partial^2\psi/\partial y^2)}{\psi}  \right\} g = 0 .\label{ma17e4} \end{equation}

We will now achieve that the expression in braces in \eqref{ma17e4} takes a form $ \sum_{j=0}^{M-2} h E_j (h) y^j - \frac{y^M}{4}$  for  power series $E_0(h), ... , E_{M-2}(h)$. Let us spell out this condition:
$$  (Q(x)+hQ_1(x,h)) + h^2 \left(\frac{dy}{dx} \right)^{2}\left[ \frac{d^2x}{dy^2} \left( \frac{dx}{dy}\right)^{-1}  \frac{(\partial\psi/\partial y)}{\psi} 
-  \frac{(\partial^2\psi/\partial y^2)}{\psi} \right] \ = \ \left(\frac{dy}{dx} \right)^{2} (\sum_{j=0}^{M-2} h E_j (h) y^j - \frac{y^2}{4}). $$


A one-page elementary calculation allows us to rewrite the $h^2$-term on the left and obtain an equation
\begin{equation} Q(x)+hQ_1(x,h) \ = \ \left( \frac{dy}{dx}\right)^2 (\sum_{j=0}^{M-2} h E_j (h) y^j - \frac{y^M}{4}) \ - \  \frac{h^2}{2}  \left[
   \frac{d^3 y}{dx^3} \left( \frac{dy}{dx}\right)^{-1} -\frac{3}{2}  \left( \frac{d^2 y}{dx^2} \right)^2 \left( \frac{dy}{dx}\right)^{-2}    \right]. \label{ma17e5} \end{equation}

Following ~\cite{AKT}, in the equation \eqref{ma17e5} we replace the independent variable $x$ by $z$, where $z=A(\int_0^x \sqrt{-Q(t)} dt)^{\frac{M}{M+2}}$. It will be convenient to choose $A=(M+2)^{\frac{2}{M+2}}$.

Quite generally, change of independent variable can be performed using formulas 
$$ \frac{dy}{dx} = \frac{dz}{dx}\frac{dy}{dz}, $$
$$  \frac{d^2 y}{dx^2} = - \frac{dy}{dz} \frac{d^2 x}{dz^2} \left( \frac{dx}{dz}\right)^{-3} + \frac{d^2 y}{dz^2}\left( \frac{dx}{dz} \right)^{-2}, $$
$$ \frac{d^3 y}{dx^3} = - 3 \frac{d^2 y}{dz^2} \frac{d^2 x}{dz^2} \left( \frac{dx}{dz} \right)^{-4} - \frac{dy}{dz} \frac{d^3 x}{dz^3} \left( \frac{dx}{dz} \right)^{-4} + 3 \frac{dy}{dz} \left( \frac{d^2 x}{dz^2}\right)^2 \left( \frac{dx}{dz}\right)^{-5} + \frac{d^3 y}{dz^3} \left( \frac{dx}{dz} \right)^{-3}. $$


The equation \eqref{ma17e5} becomes:
$$ Q(x) + hQ_1(x,h) \ = \ \left( \frac{dy}{dz}\right)^2 \left( \frac{dz}{dx}\right)^2 (\sum_{j=0}^{M-2} E_j (h) y^j - \frac{y^M}{4}) \ - \ $$ $$ \ - \  \frac{h^2}{2}  \left[
   \left( \left\{ - 3\frac{d^2 y}{dz^2} \frac{d^2 x}{dz^2} \left( \frac{d x}{dz} \right)^{-3}   
- \frac{d y}{dz} \frac{d^3 x}{dz^3} \left( \frac{d x}{dz} \right)^{-3}   
+ 3 \frac{dy}{dz} \left( \frac{d^2 x}{dz^2} \right)^2 \left( \frac{d x}{dz} \right)^{-4}  
+   \frac{d^3 y}{dz^3}\left( \frac{d x}{dz}\right)^{-2} 
 \right\} {\large /} \frac{d y}{dz}\right)\right.  $$ \begin{equation}\left.  - \frac{3}{2} \left( \left\{  - \frac{dy}{dz} \frac{d^2 x}{dz^2} \left( \frac{d x}{dz} \right)^{-2}   +   \frac{d^2 y}{dz^2}\left( \frac{d x}{dz}\right)^{-1}  \right\}  {\large /} \frac{d y}{dz}\right)^2  \right]. \label{ma17e7} \end{equation}

Multiplying both sides by $\left( \frac{dx}{dz}\right)^2$ and using
$$ \left( \frac{dx}{dz} \right)^2 = - A^{-M-2}  (\frac{M+2}{2})^2 \frac{ z^{M}}{Q(x)}  $$
we get
$$  - A^{-M-2}  (\frac{M+2}{2})^2  z^{M} + h\left( \frac{dz}{dx}\right)^{-2} Q_1(x,h) \ = \ \left( \frac{dy}{dz}\right)^2 (h \sum_{j=0}^{M-2} E_j (h) y^j - \frac{y^M}{4}) \ - \ $$ \begin{equation} \ - \  \frac{h^2}{2}  \left[
    -  \frac{d^3 x}{dz^3}  \frac{d z}{dx}    
+ \frac{3}{2}  \left( \frac{d^2 x}{dz^2} \right)^2 \left( \frac{d z}{dx} \right)^{2}  
+   \frac{d^3 y}{dz^3}\left( \frac{d y}{dz}\right)^{-1} 
   - \frac{3}{2}  \left( \frac{d^2 y}{dz^2}  \right)^2 \left( \frac{d y}{dz} \right)^{-2}  \right]. \label{ma17e8} \end{equation}
To make the structure even more transparent, put 
$$\tilde Q_1(z,h):=\sum_{n\ge 0} h^n \tilde Q_{1,n}(z) := \left( \frac{dz}{dx}\right)^{-2} Q_1(x(z),h) - h\left[ -  \frac{d^3 x}{dz^3}  \frac{d z}{dx}    
+ \frac{3}{2}  \left( \frac{d^2 x}{dz^2} \right)^2 \left( \frac{d z}{dx} \right)^{2} \right]. $$
 Then \eqref{ma17e8} becomes (recall that $A=(M+2)^{\frac{2}{M+2}}$)
\begin{equation} -  \frac{z^{M}}{4}  + h\tilde Q_1(z,h) \ = \ \left( \frac{dy}{dz}\right)^2 (\sum_{j=0}^{M-2} E_j (h) y^j - \frac{y^M}{4})   \ - \  \frac{h^2}{2}  \left[
   \frac{d^3 y}{dz^3}\left( \frac{d y}{dz}\right)^{-1} 
   - \frac{3}{2}  \left( \frac{d^2 y}{dz^2}  \right)^2 \left( \frac{d y}{dz} \right)^{-2}  \right]. \label{ma17e9} \end{equation}



From \eqref{ma17e9} we can see directly that the $h^0$ term in $y(z,h)$ equals to $z$.  Thus, we can replace the unknown $y$  of \eqref{ma17e9} by $T$ as follows: 
$$y(z,h)=z+hT(z,h). $$ 
In terms of $T(z,h)$ the equation \eqref{ma17e9} becomes (with primes denoting $\frac{d}{dz}$):
\begin{equation}  -\frac{z^M}{4} + h\tilde Q_1(z,h) \ = \ \left( 1 + hT' \right)^2 (\sum_{j=0}^{M-2} E_j (h) (z+hT)^j - \frac{(z+hT)^M}{4})   \ - \  \frac{h^2}{2}  \left[ 
   \frac{hT'''}{1+hT'}   
   - \frac{3}{2}  \frac{ ( hT''  )^2}{ \left( 1+hT' \right)^{2}}  \right], \label{ma17f1} \end{equation}
which after cancellation of $h^0$ terms becomes

$$ F(T,E_0,...,E_{M-2}) \ := \ 
  \sum_{j=0}^{M-2} E_j (h) (z+hT)^j \ - \  \frac{\sum_{\mu=1}^M \binom{M}{\mu} h^{\mu-1} z^{M-\mu}  T^\mu}{4} +     $$ \begin{equation}
+ T'\left( 2 + hT' \right) \left(- \frac{(z+hT)^M}{4} + h  \sum_{j=0}^{M-2} E_j (h) (z+hT)^j  \right)  -  \frac{h^2}{2}  \frac{T'''}{1+hT'} +  \frac{3h^3}{4}  \frac{ {T''}^2}{( 1+hT' )^{2}} - \tilde Q_1(z,h) \ = \ 0. \label{ma17f2} \end{equation}


\begin{Prop} \label{formalSolu}  The relation $$ F(T,E_0,...,E_{M-2}) =0 $$  seen as an equation on formal power series in $h$ $E_j(h)=\sum_{n\ge 0} h^n E_{j,n}$, $E_{j,n}\in \C$, and $T(x,h)=\sum_{n\ge 0} T_n(z) h^n$, has a solution with $T_n(z)$ holomorphic in any connected neighborhood of $z=0$ on which all $\tilde Q_{1,j}(z)$ are holomorphic.  
\end{Prop}

\textsc{Proof.} A simple recursion with respect to the power of $h$ can be set up similarly to  ~\cite[Th.1.1, Rmk 1.1]{AKT} once we notice that 
\begin{equation} F(T,E_0,..,E_{M-2}) = E_0 +... +E_{M-2} z^{M-2} -\frac{1}{2} z^{\frac{M}{2}} ( z^{\frac{M}{2}} T)'  - \tilde Q_1(z,h) + O(h). \label{MildlyNL} \end{equation} 
$\Box$

It is the goal of the rest of the paper to estimate the growth of term of the obtained power series in $h$.

\section{Norms on Gevrey series} \label{GevreyNorms}


The content of this section \ref{GevreyNorms} is classical.

Let $U$ be an open subset of $\C$, let $P(z,h)=\sum_{k\ge 0}  h^k p_k(z)$ be a formal power series in $h$ with holomorphic coefficients $p_k(z)\in \calO(U)$. 

We say that $P(z,h)$ is Gevrey-type on $U$ if for any compact subset $K\subset U$,  there are constants $M_K$, $\rho_K$ such that 
$$ \sup_{z\in K} |p_k(z)| \ \le M_K \rho_K^k k!. $$

\vskip2.5pc

Consider the space of those $P(z,h)=P(h)$ which do not depend on $z$, endowed with the norms which depend on a parameter $t>0$:
\begin{equation}
N_0(P,t) = \sum_{k\ge 0}  \frac{|p_k|}{k!} t^n. 
\label{NmOurs}
\end{equation}

An easy calculation establishes:

\begin{Prop} \label{BMGood}  Suppose $P(h),Q(h),$ and $t$ are such that $N_0(P,t)$, $N_0(Q,t)<\infty$. Then:\\
$$ N_0(PQ, t) \ \le \ N_0(P,t) \cdot N_0(P,t),$$
$$ N_0(P+Q,t) \ \le \ N_0(P,t) \ + \ N_0(Q,t). $$
\end{Prop}

 Let $\sfV^{t}$ be the subspace of those $P(h)$ for which $N_0(P,t)$ is finite. Analogously to the proof of completeness of $\ell^1(\C)$, one verifies that $\sfV^{t}$ is a Banach space.

Following ~\cite[p.15]{Sch}, for a Banach space $\sfA$ and  a number $\rho>0$ denote by $\sfA(\rho)$ the space of all formal series $\{ g=\sum_{j\ge 0} a_j \tau^j; \ a_j\in\sfA\}$ such that 
\begin{equation}  ||g||_\rho \ := \ \sum_{j\ge 0} ||a_j|| \rho^j \label{Sch212} \end{equation}
is finite. The space $\sfA(\rho)$ endowed with $||\cdot||_\rho$ is a Banach space (which can be proven analogously to the proof of completeness of the $\ell^1(\C)$ space). If $\sfA$ is a Banach algebra with $||ab||\le ||a||\, ||b||$ (e.g., if $\sfA=\sfV^t$), then given another such series $\tilde g=\sum_{j\ge 0} \tilde a_j \tau_j$, we have $||g\tilde g||_\rho \le ||g||_\rho\, ||\tilde g||_\rho$ for the usual product of power series.

Clearly an element of $\sfV^t(\rho)$ for $\rho,t>0$ gives rise to a Gevrey series $P(z,h)$ for $z$ in a small disc around the origin.

The following Lemma replaces the Cauchy integral formula when it comes to estimating the norm of a derivative of $\sum_j a_j \tau^j$:

\begin{Lemma} \label{ReplaceCauchyFla} If $g=\sum_{n\ge 0} a_n \tau^n \in \sfA(\rho)$ and we let $g^{(k)}:=\sum_{n\ge 0} a_n n(n-1)...(n-k+1) \tau^{n-k}$ for $k\ge 1$, then for any $\varepsilon$, $0<\varepsilon<\rho$, we have
$$ \| g^{(k)}\|_{\sfA(\rho-\varepsilon)} \le \frac{k!}{\varepsilon^k} \|g\|_{\sfA(\rho)}.$$
\end{Lemma}

\textsc{Proof.} Indeed, by the binomial formula
$$ \frac{\varepsilon^k}{k!} \sum_{n\ge 0} \|a_n\| n(n-1)...(n-k+1) (\rho-\varepsilon)^{n-k} =
 \sum_{n\ge 0} \|a_n\| \binom{n}{k} (\rho-\varepsilon)^{n-k} \varepsilon^k \le   \sum_{n\ge 0} \|a_n\| \rho^n. $$
$\Box$
 

\section{Functional analytic setup} \label{FuncAn}

Let us specify the functional spaces between which the nonlinear functional $F$ given in \eqref{ma17f2} will define a continuous map.  



\begin{Notation} \label{tau0rho0}{\rm Let us fix for the rest of the paper $\tau_0>0$, $\rho_0>0$ in such a way that $\tilde Q_1(z,h)$ defines an element of $\sfV^{\tau_0}(\rho_0)$. \\
If necessary, shrink $\rho_0$ to be $<1/M$, this will be used on p.\pageref{ma31lab5pm}. }
\end{Notation}

We define the following spaces for $0<s\le 1$, $0<t\le \tau_0$:
\begin{equation} \calX_s^t := \{ (E_0,...,E_{M-2},T) \in \underbrace{\sfV^t\times \dots \times \sfV^t}_{\text{$M-1$ factors}} \times \sfV^t(\frac{\rho_0(1+s)}{2}) \  : \ T, T',T'', T''' \in \sfV^t(\frac{\rho_0(1+s)}{2})\}  \end{equation} 
endowed with the norm \begin{equation}  ||(E_0,...,E_{M-2},T)||_{s,t} := \sum_{j=0}^{M-2} \| E_j \|_{\sfV^t} + \sum_{j=0}^3 \| \frac{d^j}{dz^j}T \|_{\sfV^t(\frac{\rho_0(1+s)}{2})};  \label{nrmXs} \end{equation}
(the reason we do not want to consider arbitrarily small neighborhoods of $z=0$ is the appearance of the radius of a neighborhood in the denominator of the estimates in Lemma \ref{LeA3zM});
\begin{equation} \calY _s^t :=  \sfV^t(\frac{\rho_0(1+s)}{2});  \end{equation} 
\begin{equation} \calZ_s^t := \{ (E_0,...,E_{M-2},T) \in \underbrace{\sfV^t\times \dots \times \sfV^t}_{\text{$M-1$ factors}} \times \sfV^t(\frac{\rho_0(1+s)}{2}) \  : \ T, T' \in \sfV^t(\frac{\rho_0(1+s)}{2})\}  \end{equation} 
endowed with the norm \begin{equation}  ||(E_0,...,E_{M-2},T)||_{s,t} := \sum_{j=0}^{M-2} \| E_j \|_{\sfV^t} + \sum_{j=0}^1 \| \frac{d^j}{dz^j}T \|_{\sfV^t(\frac{\rho_0(1+s)}{2})}.  \label{nrmZs} \end{equation}

Clearly, for every $s,t$ such that $0<s\le 1$, $0<t\le \tau_0$, $F$ defines a continuous and even analytic map  $F^t_s$  from a neighborhood $\calU_s^t$ of the origin in $\calX^t_s$ to the space $\calY^t_s$; for concreteness, $\calU_s^t=\{ (E_0,...,E_{M-1},T : t\|T'\|_{\sfV^{t}(s)}<1\}$.  

The spaces $\calZ_s^t$ are auxiliary and will be used in the proof later on.

For any fixed $t$, we have the following properties, of which 1) and 2) are trivialities and 3) follows from Lemma \ref{ReplaceCauchyFla}:

1) For $s'<s$, there are inclusions of norm $\le 1$: $\calX_s^t\hookrightarrow \calX_{s'}^t$, $\calY_s^t\hookrightarrow \calY_{s'}^t$, $\calZ_s^t \hookrightarrow \calZ_{s'}^t$. This is expressed by saying that $\calX_{(.)}^t$, $\calY_{(.)}^t$, $\calZ_{(.)}^t$ are {\it scales of Banach spaces}

2) The obvious map $\calX_s^t\stackrel{\id}{\to} \calZ_s^t$ has norm $\le 1$; 

3) For every $s'<s$, there is a map $\calZ_{s}^t \to \calX_{s'}^t$ so that
\begin{equation}  \| \calZ_{s}^t \to \calX_{s'}^t \| \le  \frac{B}{(s-s')^2}, \label{SecndDerCst} \end{equation}  where $B$ depends on $\rho_0$ but not on $t,s,s'$.





\section{The dominant term of the equation \mbox{$F(T,E_0,...,E_{M-2})=0$}} \label{secLA3}

In this section \ref{secLA3} we are studying the $h^0$ component  of the equation $F(T,E_0,...,E_{M-2})=0$, cf. \eqref{ma17f2}. The following lemma is inspired by  ~\cite[Lemma A.3]{AKT}.

\begin{Lemma} \label{LeA3zM}  Let $M\ge 1$ be an integer,  $v(z)$ a holomorphic function on $\Delta=\{ z: \, |z|<r_0\}$ with values in a Banach space $B$, and consider the following equation for an unknown holomorphic function  $u(z):\Delta\to B$ and constants $E_0,...,E_{M-2}\in B$:
\begin{equation} \left(\frac{ z^M}{2} \frac{d}{dz} + \frac{M}{4}z^{M-1} \right) u(z) \ = \  E_0+zE_1+...+z^{M-2}E_{M-2} +v(z). \label{zM19} \end{equation}
 Then \eqref{zM19} has a unique solution $(u(z),E_0,...,E_{M-2})$ and moreover for any $r$, $0<r<r_0$: 
\begin{equation} ||E_j||_B \ \le \ \frac{1}{r^j} || v ||_{B(r)}, \ \ j=0,...,M-2;
\label{mA20} \end{equation}
\begin{equation} ||u||_{B(r)} \ \le \ \frac{4}{r^{M-1}} ||v||_{B(r)};
\label{mA21} \end{equation}
\begin{equation} \left| \left| \frac{du}{dz}\right|\right|_{B(r)} \ \le \ \frac{2}{r^M } ||v||_{B(r)} .
\label{mA22} \end{equation}
\end{Lemma}


\textsc{Proof}. 
%
Rewrite the equation \eqref{zM19} as follows:
\begin{equation} \frac{1}{2} z^{\frac{M}{2}} (z^{\frac{M}{2}}u(z))' = E_0+zE_1+...+z^{M-2}E_{M-2} +v(z). \label{z701} \end{equation}
If $v(z)=v_0+v_1 z + v_2 z^2 +...$, then \eqref{z701} can be rewritten as
\begin{equation}  \frac{1}{2}  (z^{\frac{M}{2}}u(z) )'= \sum_{j=0}^{M-2} z^{j-\frac{M}{2}} (E_j + v_j) + \sum_{j=M-1}^\infty z^{j-\frac{M}{2}} v_j . \label{z703} \end{equation}
If $u(z)$ is to be holomorphic, the RHS of \eqref{z703} should not have any terms with $z^{\le -1}$, hence 
$$ v_j = -E_j, \ \ \ \text{if} \ j-\frac{M}{2} \le -1. $$
Inserting this, we have
\begin{equation}  \frac{1}{2}  (z^{\frac{M}{2}}u(z) )'= \sum_{j\in \N_0; \,  \frac{M}{2}-1\le j \le M-2} z^{j-\frac{M}{2}} (E_j + v_j) + \sum_{j=M-1}^\infty z^{j-\frac{M}{2}} v_j , \label{z705} \end{equation}
or
\begin{equation}   u(z)= \sum_{j\in \N_0; \,  \frac{M}{2}-1\le j \le M-2} \frac{2}{j-\frac{M}{2}+1}z^{j-M+1} (E_j + v_j) + \sum_{j=M-1}^\infty \frac{2}{j-\frac{M}{2}+1} z^{j-M+1} v_j.  \label{z715} \end{equation}
We conclude that the solution is necessarily of the form
\begin{equation} E_j = -v_j, \ \ \ j=0,...,M-2, \label{ma274pm} \end{equation}
$$ u(z) = \sum_{j=M-1}^\infty \frac{2}{j-\frac{M}{2}+1} v_j z^{j-M+1}. $$

As $j-\frac{M}{2}+1\ge \frac{1}{2}$ for $j\ge M-1$, we have
$$ \| u(z) \|_{B(r)} \le 4 \sum_{j=M-1}^\infty \| v_j\|_B r^{j-M+1} \le 4 r^{1-M} \sum_{j=M-1}^\infty \|v_j\|_B r^j \le 4 r^{1-M} \| v(z)\|_{B(r)} .$$

Further,
$$ u'(z) = 2 \sum_{j=M-1}^\infty \frac{j-M+1}{j-\frac{M}{2}+1} v_j z^{j-M}, $$
hence, as the fraction is always $\le 1$,
$$ \| u'(z) \|_{B(r)} \le 2 \sum_{j=M}^\infty \| v_j\|_B r^{j-M} \le 2 r^{-M} \| v(z)\|_{B(r)}. $$ 

From \eqref{ma274pm} we obviously get \eqref{mA20} . $\Box$

\section{Newton's method} \label{NMblah}

This section \ref{NMblah} contains no precise results, it will not be referred to except for purposes of intuition. 

\subsection{Newton's method in Banach spaces.} \label{NMsec}
 
Newton's method of tangents for solving scalar nonlinear equations appears in almost all elementary calculus textbooks. In this section \ref{NMsec} we review the Newton's method in Banach spaces following  ~\cite[Ch. XVIII]{Kanto}.

Suppose $\calX$, $\calY$ are two Banach spaces, $x^{(0)}\in \calX$, $\calU\subset \calX$ is an open neighborhood, $\calF: \calU\to \calY$ a continuous map admitting two continuous and bounded Fr\'echet derivatives.   Suppose that for all $x\in \calU$ the inverse of the Fr\'echet derivative, $(dF_x )^{-1}:\calY \to \calX$ is uniformly bounded, i.e. $\|(dF_x)^{-1}\|\le A$.  

By Newton's method we mean the following iterative procedure of finding a zero of $\calF$:\\
1) $x^{(0)}$ is given; $y^{(0)}=F(x^{(0)})$;  \\
2) for $k\ge 0$, define $w^{(k)} = -(dF_{x^{(k)}})^{-1}(y^{(k)})$; \\
3) $x^{(k+1)} = x^{(k)} + w^{(k)}$; assume or prove that $x^{(k+1)}\in \calU$;  $y^{(k+1)} = F(x^{(k+1)})$. 

Assume that $x^*$ is such that $F(x^*)=0$. Then 
$$ y^{(k+1)} = F(x^{(k+1)}) = \underbrace{F(x^{(k)}) + [dF_{(x^{k})}] (w^{(k)})}_{=0} + O(\| w^{(k)}\|^2)  = O (\| y^{(k)} \|^2 ) $$
which is smaller than $\| y^{(k)}\|$ if $y^{(k)}$ was small already. So the method converges if $x^{(0)}$ was sufficiently close to $x^*$ and, loosely, far enough from the boundary of $\calU$. Formalization of these conditions can be found in ~\cite[Ch.XVIII, {\S}1.5]{Kanto}. 

 

\subsection{Newton's method leads to Gevrey expansions} \label{dF0}

In this subsection \ref{dF0} is to motivate that factorially divergent or, for brevity, Gevrey expansions naturally arise when we attempt to apply Newton's method for solving our equation $F(E_0,...,E_{M-2},T)=0$.


The Newton's method requires inverting the Fr\'echet derivative of $F$. Let us for simplicity consider the case $M=2$ and try to invert the Fr\'echet derivative of $F$ at the point $(E_0,T)=0$. 
If  $(\calE_0,\calT)$ is a tangent vector, then 
$$ dF_0 (\calE_0, \calT) = \calE_0 - \frac{z^2}{2}\calT' - \frac{z}{2}\calT - \frac{1}{2} h^2 \calT'''. $$
Then $dF_0$ as an operator $\calX^t_s\to \calY^t_s$ can be written as $G+h^2 H$, where $H$ is bounded and $$ G \ : \ (\calE_0,\calT)\mapsto \calE_0-\frac{z^2}{2}T' - \frac{z}{2}\calT \ : \ \calX^t_s\to \calY^t_s$$ is inverted by an operator $K:\calY_s \to \calX_{s'}$, $s'<s$, of the norm $\le \frac{c}{(s-s')^2}$, cf. lemma \ref{LeA3zM} and \eqref{SecndDerCst}.   
%
%
Let us attempt to write the inverse of $dF_0$ as
$$ (dF_0)^{-1} = \sum_{n\ge 0} h^{2n} (-1)^n (KH)^n K $$
where $(KH)^n K$ is 
\begin{equation} \calY^t_s \stackrel{K}{\to} \calX^t_{s_1} \stackrel{H}{\to} \calY^t_{s_1} \stackrel{K}{\to} \calX^t_{s_2} \stackrel{H}{\to} ... \stackrel{H}{\to} \calY^t_{s_n} \stackrel{K}{\to} \calY^t_{s'} \end{equation}
and where we can think of $s_n<...<s_1$ as arbitrary numbers subject to $s'<s_n$ and $s_1<s$. Choosing $s_1,..,s_n$ to maximize the product for $(s-s_1)(s_1-s_2)...(s_n-s')$ for fixed $s,s'$, we can prove an estimate 
$$ \| (KH)^n K \ : \ \calX^t_{s} \to \calY^t_{s'} \| \le c'' \frac{(2n)! (c')^{2n}}{(1-s)^{2n}} $$
for some new constants $c'$, $c''$, but no dramatically better estimates are available.
Thus, given $v\in \calY^t_1$, we can hope to represent $(dF_0)^{-1}v \in \calX^t_s$ by the following expansion whose convergence needs to be discussed separately:
\begin{equation} (dF_0)^{-1} v = \sum_{n\ge 0} h^n u_{n,s}, \ \ \ \| u_{n,s}\|_{\calX^t_s} \le C' \frac{n! C^n}{(1-s)^n}, \label{dF0invV} \end{equation}
with constants $C',C$ independent of $s$. This shows that expansions of type \eqref{dF0invV} appear naturally in our problem and motivates our decision to formulate intermediate results in terms of them.


\section{Calculus of Gevrey expansions}  \label{Gcalc}

\subsection{Definitions} \label{GcalcDef}

By a Gevrey expansion we mean, vaguely, a formal expansion 
$ \sum_n \rmh^n u_n$, where the straight-font $\rmh$ is a formal variable, and $u_n$ are elements of some Banach space such that $\| u\|\le c_0 c^n n!$ for some constants $c_0$ and $c^n$, compare  \eqref{dF0invV}. We choose however to avoid this notion in the mathematically precise body of the article. Instead, the following seems to be a more convenient terminology.

In this section \ref{Gcalc}, let $\calX_s$, $\calY_s$, $0<s\le 1$ be two arbitrary scales of Banach spaces, i.e. we suppose that for $s'<s$ there is an inclusion $\calX_s\hookrightarrow \calX_{s'}$, $\calY_s\hookrightarrow \calY_{s'}$ of norm $\le 1$. In this section \ref{Gcalc}, $\calX_s, \calY_s$ have a priori nothing to do with $\calX^t_s$, $\calY^t_s$ in the rest of the paper.

Suppose that for every $s$ there is a map $F_s: \calX_s \supset \calU_s \to \calY_s$. We will say that the family $F=(F_s)$ is {\it compatible with inclusions} if:\\ 
a) for $s'<s$, the map $\calX_{s}\hookrightarrow \calX_{s'}$ maps $\calU_s$ into $\calU_{s'}$; \\
b) for $s'<s$, $(\calY_{s}\hookrightarrow \calY_{s'})\circ F_s = F_{s'} \circ (\calX_{s}\hookrightarrow \calX_{s'})$.  \\
For a family of linear maps $F_s$ we always think of $\calU_s$ as equal to $\calX_s$.

Analogously we define what it means for a family of maps $(G_{s,s'}:\calX_s \to \calY_{s'})_{s'<s}$ to be compatible with inclusions.

We denote by $\mathring\calX_1$ the inductive limit of $\calX_s$ for $s<1$:
$$ \mathring\calX_1 := \{ (g_s)_{0<s<1} \ : \ g_s \in \calX_s, \ \text{and for all $s<s'$} \ (\calX_{s'}\hookrightarrow \calX_s)(g_{s'})=g_s\}. $$

We assume that on every $\calX_s$ there is an action of a linear operator $h:\calX_s\to \calX_s$, compatible with inclusions, satisfying $\| h^n: \calX_s\to \calX_s\|\le \frac{\tau}{n!}$, and similar for $\calY_s$, for some fixed number $\tau>0$ independent of $s$. We will abuse notation by writing $|h|$ instead of $\tau$. 

We call a linear operator $G_s:\calX_s\to \calY_s$ (or a family of such operators) {\it $h$-linear} if $G_s \circ h = h\circ G_s$. 

We finish this subsection by stating the following combinatorial inequality:

\begin{Lemma} \label{AKTLeA4} {\rm ~\cite[Lemma A.4]{AKT}} The following inequality holds for all positive integers $j,k$ satisfying $k\le j$:
$$ \sum_{\footnotesize \begin{array}{c} j_1+...+j_k=j \\ j_1,\dots,j_k\ge 1\end{array}} j_1! j_2! \dots j_k! \le 4^{k-1}(j-k+1)!. $$
\end{Lemma} 


\subsection{Function evaluated on a Gevrey expansion gives a Gevrey expansion}




%

For $k\in\Z_{\ge 0}$, we call $\Omega_k:\calX_s\to \calY_s$ an $h$-$k$-linear map of norm $\le 1$ if $\Omega_k (v) = \tilde \Omega_k(v,v,...,v)$ where $\tilde \Omega_k: \calX_s\times \calX_s\times \dots \times \calX_s \to \calY_s $ is symmetric, satisfies $||\tilde \Omega(v_1,...,v_k)||_{\calY_s}\le ||v_1||_{\calX_s}\, ||v_2||_{\calX_s}\, ... \, ||v_k||_{\calX_s}$ and, moreover, $\tilde \Omega_k(...,v_{j-1},hv_j,v_{j+1},...) =h\tilde\Omega_k(...,v_{j-1},v_j,v_{j+1},...)$.

\begin{Lemma} \label{InsertL} Let $f: \calX_s\to \calY_s$ be a map compatible with inclusions, $x_0\in \calX_1$; $\alpha\in \R_{>0}$. Suppose that for any $w \in \calX_s$ such that  $||w||_s < \frac{1}{\alpha} $ 
\begin{equation} \calY_s \ni f(x_0+w) = f(x_0) + \sum_{j\ge 1} \alpha^j \Omega_j (w), \label{fe12d} \end{equation}
where $\Omega_{j,s}:\calX_s\to \calY_s$ is an $h$-$j$-linear map of norm $\le 1$, compatible with inclusions. \\
Let for all $n\ge 1$, $g_n = (g_{n,s})\in \mathring\calX_1$ are such that
\begin{equation} ||g_{n,s}||_{\calX_s} \le C_0 \frac{C^n n!}{(1-s)^n}. \label{ma6a} \end{equation} 
Then there are $H_{n}=(H_{n,s})\in \mathring \calY_1$, $n\ge 1$ satisfying
$$ ||H_{n,s}||_{\calY_s} \le \beta(C_0\alpha)\frac{C^n n!}{(1-s)^n} $$ 
for $\beta(t)=te^{4t}$,  such that if 
\begin{equation} \frac{ |h|C}{1-s} <1  \ \ \ \text{and} \ \ \ \sum_{n\ge 1} \frac{|h|^n C^n}{(1-s)^n} < \frac{1}{C_0\alpha}, \label{con11f} \end{equation}
then the series  \begin{equation} W_s = \sum_{n\ge 1} h^n g_{n,s} \label{fe12b} \end{equation}
converges in $\calX_s$  and  \begin{equation} f(x_0+W_s) - f(x_0) = \sum_{n\ge 1} h^n  H_{n,s}.  \label{fe12c} \end{equation}
\end{Lemma}

\textsc{Proof.}  Motivated by a ``formula" 
{ \begin{equation} \text{ `` } \ \  f(x_0+W_s) -f(x_0) = \sum_{k\ge 1} a_k  \tilde \Omega_k \left( \sum_{n_1\ge 1} h^{n_1} g_{n_1,s}, ..., \sum_{n_k\ge 1} h^{n_k} g_{n_k,s} \right)  \ \ \ \  \ \ \ \ \ \ $$
$$ \ \ \ \ \ \ \ \ \  \ = \ \sum_{k\ge 1} a_k \sum_{m\ge 1} h^m\sum_{\footnotesize \begin{array}{c} n_1+...+n_k=m \\ n_1,...,n_k\ge 1 \end{array}} \tilde \Omega_k \left( g_{n_1,s}, ..., g_{n_k,s} \right) \ \text{"} \label{fe12a} \end{equation} }
to which we will give an analytic meaning in a moment, consider 
$$H_{m,s} := \sum_{k\ge 1} a_k \sum_{\footnotesize \begin{array}{c} n_1+...+n_k=m \\ n_1,...,n_k\ge 1 \end{array}} \tilde \Omega_k \left( g_{n_1,s}, ..., g_{n_k,s} \right).$$
We have, with help of lemma \ref{AKTLeA4}, 
$$ ||H_{m,s}||_{\calY_s} \le \sum_{k\ge 1} (C_0\alpha)^k \frac{C^m}{(1-s)^m} \sum_{\footnotesize \begin{array}{c} n_1+...+n_k=m \\ n_1,...,n_k\ge 1 \end{array}} n_1!...n_k! \le $$
$$ \le \frac{C^m}{(1-s)^m} \sum_{k\ge 1} (C_0\alpha)^k  4^{k-1} (m-k+1)!  \le (C_0\alpha) \frac{C^m}{m! (1-s)^m} \sum_{k\ge 0} (C_0 \alpha)^{k}  4^{k} \frac{1}{k!} \le (C_0 \alpha) e^{4C_0 \alpha} \frac{ C^m m!}{ (1-s)^m},    $$
compare ~\cite[(A.50)]{AKT}.

Then, assuming \eqref{con11f} and hence $\left\| \frac{h^n C^n}{(1-s)^n} \right\| \le \frac{1}{n!}$, 
$\sum_{m\ge 1} h^m H_{m,s}$  as well as \eqref{fe12b}  are absolutely convergent, so $f(x_0+W_s)$ can be written in terms of Taylor series \eqref{fe12d}, and hence \eqref{fe12c} holds by the algebraic manipulations of \eqref{fe12a}.
$\Box$


\subsection{Increment of a function is well approximated by the first derivative}

Let us now give a meaning to the formula $f(x+h^b v)= f(x)+[f'(x)](h^b v) + h^{2n} O(v^2)$ for $x=x_0+g$ in the context of Gevrey expansions. 


\vskip2.5pc


Let $\beta_1(t) = \sum_{k\ge 0} \frac{ (k+1) 4^k}{k!} t^{k+1}$, $\beta_2(t) =\sum_{k\ge 1} \frac{t^{k+1} 4^k}{k!}$.

\begin{Lemma} \label{DiffWellDef}  Suppose $f:\calX_s\to \calY_s$, $x_0$, $\alpha$,$\Omega_j$ are as in \eqref{fe12d}, and $g_n\in \mathring \calX_1$, $n\ge 1$, satisfy \eqref{ma6a}.  Assume $C\ge 1$. \\
Then for any integer $b\ge 0$ and any sequence $v_k\in \mathring \calX_1$ 
with $\| v_{k,s}\|_{\calX_s} \le A \frac{k! C^k}{(1-s)^k}$, there are : \\ 
a) a sequence of elements $w_j\in \mathring \calY_1$ with $\|w_{j,s}\| \le \frac{A}{C_0} \beta_1(C_0\alpha) \frac{j! C^j}{(1-s)^j}$ satisfying the property:\\
under condition 
\begin{equation} \frac{|h|C}{(1-s)}<1, \ \ \ \sum_{n=1}^\infty \frac{h^n C^n}{(1-s)^n} < \frac{1}{C_0 \alpha}, \label{ma8b} \end{equation} we have two vector defined by an absolutely convergent sums
$$ v_s :=h^b \sum_{k=1}^\infty h^{k} v_{k,s} \in \calX_s,  \ \ g_s:=\sum_{k\ge 1} h^k g_{k,s} \in \calX_s$$
and such that
\begin{equation}  f'(x_0+g_s)(v_s) = \sum_{j\ge 1} j \alpha^{j} \tilde \Omega_j(g_s,...,g_s,v_s)  \ = \ \sum_j h^j w_{j,s}  \ \ \  \text{in} \ \calY_s; \label{ma8a} \end{equation}
b) a sequence of elements $u_j\in \mathring\calY_1$, $j\ge 1$, with 
$$ \| u_{j,s} \|_{\calY_s} \le  \beta_2(\alpha(C_0+A)) \frac{j! C^j}{(1-s)^j} $$ 
satisfying the following property: \\
under conditions \begin{equation} \frac{|h|C}{(1-s)}<1, \ \ \ \sum_{n=1}^\infty \frac{h^n C^n}{(1-s)^n} < \frac{1}{ (C_0+A)\alpha} \label{ma8c} \end{equation} 
we have
\begin{equation} f(x_0+g_s+v_s)-f(x_0+g_s)-[f'(x_0+g_s)](v_s) \ = \  h^{2b} \sum_{j\ge 1} h^j u_{j,s} \ \ \ \text{in} \ \calY_s. \label{mar6o} \end{equation}
\end{Lemma}

\textsc{Proof.} 

\underline{Preliminary remark.}  If $v$ is a formal in $h$ written as  $v=h^b \sum_{n\ge 1} h^k v_k$ as in the statement of the lemma, then it can also be written $v=\sum_{m\ge b+1} h^m \tilde v_m$ with $\tilde v_m = v_{m-b}$ and hence satisfying $\| \tilde v_{m,s} \|\le A\frac{C^{m-b} (m-b)!}{(1-s)^{m-b}}  \le  A\frac{C^{m} m!}{(1-s)^{m}} $  if $C\ge 1$. 

\underline{Proof of a).}
Rewrite \eqref{ma8a} as the following ``formula" whose analytic meaning will be clarified in a moment: \begin{equation} \text{``} \ \sum_{j\ge 1} \sum_{k_1,..,k_{j-1},k_j\ge 1} j \alpha^{j}  h^{k_1+...+k_{j-1}} \tOmega_j(g_{k_1},...,g_{k_{j-1}},\tilde v_{k_j}) \ = \  \sum_{n\ge 1} h^n w_n .\ \text{"} \label{mar6d} \end{equation}
Thus we define
$$ w_n = \sum_{j\ge 1} \sum_{\footnotesize \begin{array}{c} k_1,..,k_{j-1},k_j\ge 1 \\ k_1+...+k_{j-1}+k_j=n \end{array} } j C_0^{j-1} A \alpha^j \Omega_j(g_{k_1},...,g_{k_{j-1}}, \tilde v_{k_j} ), $$
and with help of Lemma \ref{AKTLeA4}
$$ \| w_{n,s} \|_{\calY_s}   \le  \sum_{j\ge 1} \sum_{\footnotesize \begin{array}{c} k_1,..,k_{j}\ge 1 \\ k_1+...+k_{j}=n \end{array} } j \alpha^j C_0^{j-1} A  \frac{C^n k_1!...k_{j}! }{(1-s)^{n}}  \le $$ 
$$ \le  \sum_{j\ge 1}  j \alpha^j  C_0^{j-1} A \frac{C^n 4^{j-1} (n-j+1)! }{(1-s)^{n}}  \le \frac{A}{C_0} \ \left\{ \sum_{j'\ge 0}  (j'+1) 4^{j'} (C_0 \alpha)^{j'+1} \frac{1}{j'!}  \right\} \frac{C^n  n! }{(1-s)^{n}}.  $$ 

Under condition \eqref{ma8b}, the equation \eqref{mar6d} is literally true in $\calY_s$, both sides being absolutely convergent series in which we can change the order of summation, hence \eqref{ma8a}. 


\underline{Proof of b)}. Let us now formally write 
$$ f(x_0+g+v) - f(x_0+g) - [f'(x_0+g)](v) \ \text{``}=\text{"} \ \ \ \ \ \ $$
$$ \ = \ \sum_{j=1}^\infty \alpha^j \left[ \Omega_j(g+v) - \Omega_j(g) - j\tilde \Omega_j(g,...,g,v) \right] \ = \ $$
(use symmetry and polylinearity of $\tilde\Omega_j$) 
$$ \ = \ \sum_{j=2}^\infty \alpha^j\sum_{k_1,...,k_j\ge 1} h^{k_1+...+k_j} \left[ \sum_{\sigma=2}^j \binom{j}{\sigma} \tilde \Omega_j(g_{k_1}, ... ,g_{k_\sigma}, \tilde v_{k_{\sigma+1}},...,\tilde v_{k_j}) \right] \ =$$
(since $\tilde v_{k}=0$ for $k=0,...,b$ )
$$ \ = \ \sum_{j=2}^\infty \alpha^j\sum_{\footnotesize \begin{array}{c} k_1,...,k_{j-2}\ge 1; \\ k_{j-1},k_j \ge b+1 \end{array}} h^{k_1+...+k_j} \left[ \sum_{\sigma=2}^j \binom{j}{\sigma} \tilde \Omega_j(g_{k_1}, ... ,g_{k_\sigma}, \tilde v_{k_{\sigma+1}},...,\tilde v_{k_j}) \right] \ = $$
$$ \ = \ \sum_{j=2}^\infty \alpha^j\sum_{\footnotesize \begin{array}{c} k_1,...,k_{j-2}\ge 1; \\ k_{j-1},k_j \ge b+1 \end{array}} h^{k_1+...+k_j} \left[ \sum_{\sigma=2}^j \binom{j}{\sigma} \tilde \Omega_j(g_{k_1}, ... ,g_{k_\sigma}, \tilde v_{k_{\sigma+1}},...,\tilde v_{k_{j-2}}, v_{k_{j-1}-b}, v_{k_j-b}) \right] \ = $$
$$ \ = \ \sum_{j=2}^\infty \alpha^j\sum_{ k_1,...,k_{j}\ge 1} h^{k_1+...+k_j+2b} \left[ \sum_{\sigma=2}^j \binom{j}{\sigma} \tilde \Omega_j(g_{k_1}, ... ,g_{k_\sigma}, \tilde v_{k_{\sigma+1}},...,\tilde v_{k_{j-2}}, v_{k_{j-1}}, v_{k_j}) \right]. $$

This motivates the definition: 
$$ u_n = \sum_{j=2}^\infty \alpha^j\sum_{ \footnotesize \begin{array}{c} k_1,...,k_{j}\ge 1 \\ k_1+...+k_j=n \end{array} }  \left[ \sum_{\sigma=2}^j \binom{j}{\sigma} \tilde \Omega_j(g_{k_1}, ... ,g_{k_\sigma}, \tilde v_{k_{\sigma+1}},...,\tilde v_{k_{j-2}}, v_{k_{j-1}}, v_{k_j}) \right] . $$
We have:
$$ \| u_{n,s} \|_{\calY_s} \le \sum_{j=2}^\infty \alpha^j\sum_{ \footnotesize \begin{array}{c} k_1,...,k_{j}\ge 1 \\ k_1+...+k_j=n \end{array} }  \left[ \sum_{\sigma=2}^j \binom{j}{\sigma} \frac{C_0^\sigma A^{j-\sigma} C^n}{(1-s)^{n} } {k_1}! ...{k_j}! \right] \le $$
$$  \le \sum_{j=2}^\infty \alpha^j (C_0+A)^j 4^{j-1} (n-j+1)!    \frac{ C^n}{(1-s)^{n} } \le \left\{ \sum_{j=2}^\infty \alpha^j (C_0+A)^j 4^{j-1} \frac{1}{(j-1)!} \right\}    \frac{ C^n n!}{(1-s)^{n} }.$$
We conclude that \eqref{mar6o} holds under condition \eqref{ma8c} in the same way as we did in part a). $\Box$


\subsection{Inverse of a regularly perturbed linear operator} \label{dGminus1}

In this subsection \ref{dGminus1} we will generalize the familiar statement that if an operator $\bfA$ has a bounded inverse, then its small perturbation $\bfA+\bfB$ has a bounded inverse $\bfA^{-1}(1+\bfB \bfA^{-1} + \bfB \bfA^{-1} \bfB \bfA^{-1} +...$ if $\|\bfB\| \| \bfA^{-1}\|<1$. We show that a similar statement holds if $\bfB$ is replaced by a factorially divergent series in $h$ with operator coefficients.  



\begin{Lemma} \label{BenignInv} Suppose given a sequence of $h$-linear compatible with inclusions operators: 
\begin{equation} \calH_{n,s} : \calX_s \to \calY_s \ : \  ||\calH_{n,s}||\le C_1 \frac{n! C^n}{(1-s)^n}, \ \ \ n\ge 1, \ 0<s< 1. \label{as12ma} \end{equation}
Suppose that we also  have $h$-linear compatible with inclusions operators $\calF_s:\calX_s\to\calY_s$, with inverse $\calG_s:\calY_s\to \calX_s$, $||\calG_s||\le A$, which is also  $h$-linear and compatible with inclusions. \\ 
Then there exist operators, $h$-linear and compatible with inclusions,
$$ \calL_{n,s} : \calY_s\to \calX_s \ : \ ||\calL_{n,s}:\calY_s\to\calX_s || \le \max\{ A, (C_1 A)^2 e^{4(C_1 A)}\} \frac{C^n n!}{(1-s)^n} , \ \ \ n\ge 0, \ \ 0<s< 1$$
satisfying the following property: \\
with notation $\calH_0:=\calF$, 
\begin{equation} \sum_{m,j\ge 0; \ j+m=N} \calH_{j,s} \calL_{m,s} = \left\{ \begin{array}{ll} 0, &  \text{if} \ N\ge 1; \\ \id, & \text{if} \ N=0. \end{array} \right. \label{ma20e1} \end{equation} 
 %
\end{Lemma}

%
%

\textsc{Proof.} 
Motivated by the ``formulas" 
\begin{equation} \text{``} \ (\calF +\calH)^{-1} = ((1+ \calH\calF^{-1}) \circ \calF) ^{-1} = \calG \circ (1+\calH \calG)^{-1} = \sum_{k\ge 0} (-1)^j \calG (\calH \calG)^k \ \text{"} \label{ma11s1} \end{equation}
and,  for $k\ge 1$ 
$$ \text{``} \ ( \calH \calG)^k = \sum_{m\ge 1}  h^m \sum_{\footnotesize \begin{array}{c} n_1+...+n_k = m \\ n_1,...,n_k\ge 1 \end{array}} \calH_{n_1,s} \calG \calH_{n_2,s}\calG ... \calH_{n_k,s}\calG\ \text{"},  $$
whose analytic meaning will be clarified shortly, we 
put $$ \calK_{k,m,s} := \sum_{\footnotesize \begin{array}{c} n_1+...+n_k = m \\ n_1,...,n_k\ge 1 \end{array}} \calH_{n_1,s} \calG \calH_{n_2,s}\calG ... \calH_{n_k}\calG. $$
We have 
$$ ||\calK_{k,m,s} || \le \sum_{\footnotesize \begin{array}{c} n_1+...+n_k = m \\ n_1,...,n_k\ge 1 \end{array}} \frac{C_1 n_1 ! C^{n_1}}{(1-s)^{n_1}} A  \frac{C_1 n_2 ! C^{n_2}}{(1-s)^{n_2}} A  ... \frac{C_1 n_k ! C^{n_k}}{(1-s)^{n_k}} A \le  $$
$$ \ \ \ \ \ \  \le (C_1 A)^{k} \frac{C^m}{(1-s)^m} \sum_{\footnotesize \begin{array}{c} n_1+...+n_k = m \\ n_1,...,n_k\ge 1 \end{array}} n_1!...n_k!  \le  $$
(use lemma \ref{AKTLeA4})
$$ \ \ \ \ \ \  \le (C_1 A)^{k} \frac{C^m}{(1-s)^m} 4^{k-1} (m-k+1)! .  $$

Put for $m\ge 1$:
\begin{equation} \calL_{m,s} := \sum_{k\ge 1} (-1)^{k} \calG \calK_{k,m,s}; \label{ma19ni} \end{equation}
then
$$ ||\calL_{m,s} || \le (C_1 A) \sum_{k\ge 1} (C_1 A)^{k} \frac{C^m}{(1-s)^m} 4^{k-1} (m-k+1)! \le \ \ \ \ \ $$ $$ \ \ \ \ \le \ (C_1 A)^2  \frac{C^m}{(1-s)^m}  \sum_{k'\ge 0} (C_1 A)^{k'} 4^{k'} (m-k')! \le (C_1 A)^2  \frac{C^m m!}{(1-s)^m  } e^{4(C_1 A)}. $$
Put further $\calL_{0,s}=\calG$. 

One proves \eqref{ma20e1} by writing out the LHS in terms of $\calH_j$ and $\calG$ and canceling; absolute convergence in \eqref{ma19ni} justifies manipulations with infinite sums. 
  $\Box$

\normalsize


\subsection{The inverse of a singularly perturbed linear operator}




We will use the following immediate corollary of ~\cite[(A.31)]{AKT}
\begin{equation} \sum_{j_1+j_2=j; \ j_1,j_2\ge 0} j_1! j_2! \le 2j! + 4(j-1)! \le 6 j! , \ \ \ \text{if} \ j\ge 0. \label{fe15in6} \end{equation}





{\bf Remark.} Let us intuitively motivate the setup of the lemma \ref{EvilInv} below. We start with an  operator-valued Gevrey expansion $\calF=\sum_{n\ge 0} \rmh^n \calF_{n}$ which has a one-sided inverse in the form of the operator-valued Gevrey expansion $\calG=\sum_{n \ge 0} \rmh^n \calG_n$. We consider a third operator-valued Gevrey expansion $\calH=\sum_{n\ge 0} \rmh^n \calH_{n}$. We are solving a linear equation $(\calF + \rmh^a \calH) u = v$, where $v$ is a known vector-valued Gevrey expansion and $u$ is the Gevrey expansion to be found. We speak, in the title of this subsection, of a singularly perturbed linear operator because  $\calG$ does not act between $\calY_s$ and $\calX_s$ with the same $s$.

\vskip2.5pc

\begin{Lemma} \label{EvilInv}   Let $a\in\{1,2,..\}$, $C''\ge 1$, and $C\ge \max\{ 1, 36a^{a-2} BC'' e^a\}$. \\ Suppose  we are given  compatible with inclusions $h$-linear operators\begin{equation} \calF_{n,s}:\calX_s\to\calY_s \ \ : \ \ || \calF_{n,s}||_{\calX_s,\calY_s} \le C' \frac{n! C^n}{(1-s)^n}, \ \ \ 0<s<1, \ \ n\ge 0 \label{ma125pm} \end{equation}
 and   compatible with inclusions $h$-linear operators
$$ \text{for} \ n\ge 0, \ \ \text{for} \ s'>s, \ \ \  \calG_{n,s',s} \ : \calY_{s'}\to \calX_s  \ \ \ \  \| \calG_{n,s',s}\|\le \frac{B}{(s'-s)^a} \frac{n! C^n}{(1-s')^n}.$$ 
so that
\begin{equation} \sum_{m,j\ge 0; \ j+m=N} \calF_{j,s} \calG_{m,s',s} = \left\{ \begin{array}{ll} 0, &  \text{if} \ N\ge 1; \\ \calY_{s'} \stackrel{\id}{\hookrightarrow} \calY_s, & \text{if} \ N=0. \end{array} \right. \label{ma20e2} \end{equation} 
 Suppose given another family of $h$-linear operators compatible inclusions
$$ \calH_{n,s} : \calX_{s} \to  \calY_s  \ \ \ \ ||\calH_{n,s}||_{\calX_{s}, \calY_s}\le  C''\frac{n! C^n}{(1-s)^n}   $$  
and elements $v_{n} = (v_{n,s})\in \mathring \calX_1$, $n\ge 0$, satisfying $$|| v_{n,s}||\le \frac{n! C^n}{(1-s)^n}.$$

Then there are $u_n= (u_{n,s})\in \mathring\calX_1$ such that 
$$ || u_{n,s} ||_{\calX_s}  \le \frac{C^{n+a+1} (n+a+1)!}{(1-s)^{n+a+1}}, \ \ \ n\ge 0 $$
so that: \\ 
If 
\begin{equation}  \frac{|h| C}{1-s} <1  \label{EvilConerg} \end{equation}
then $\calF_{s} =\left( \sum_{n\ge 0} h^n \calF_{n,s} \right) ,  \calH_{s} =\left( \sum_{n\ge 0} h^n \calH_{n,s} \right) : \calX_s\to \calY_s$ are bounded operators, $v=\sum_{n\ge 0} h^n v_{n,s}\in\calY_s$ and $u=\sum_{n\ge 0} h^n u_{n,s} \in \calX_s$ are absolutely convergent series, and 
\begin{equation} (\calF + h^a \calH) u = v  \ \ \ \text{in} \ \calY_s. \label{EvilConclusion} \end{equation}
\end{Lemma}

\textsc{Proof.} \underline{Step 1.} (Algebraic construction of the solution) 

Guided by the ``formula"
$$\text{``} \ \  (\calF+h^a \calH)^{-1} = [(1+h^a \calH \calF^{-1})\calF]^{-1} = \calF^{-1}  (1 -  h^a\calH \calF^{-1} + h^{2a}(\calH \calF^{-1})^2 + ... ) \ \ \text{"} $$
which however does not make literal sense, we take 
%
\begin{equation}\calG (\calH \calG)^n v \text{``}=\text{"} \sum_{\footnotesize \begin{array}{c} k_0,...,k_n\ge 0 \\   \ell_1,...,\ell_n,m\ge 0\end{array}} h^{k_0+...+k_n+\ell_1+...+\ell_n+m} \calG_{k_0,s_1,s_0} \calH_{\ell_1,s_1} \calG_{k_1,s_2,s_1} ... \calH_{\ell_n,s_n} \calG_{k_n,s_{n+1},s_n} v_{m,s_{n+1}} \label{ap24c} \end{equation}
with $s_0=s$ and we reserve the right to choose other $s_j$'s differently in each summand subject only to the condition $s_0<s_1<...<s_{n+1}<1$ (compositions will not depend on these choices).
%
%
We have 
$$ || \calG_{k_0,s_1,s_0} \calH_{\ell_1,s_1} \calG_{k_1,s_2,s_1} ... \calH_{\ell_n,s_n} \calG_{k_n,s_{n+1},s_n} v_{m,s_{n+1}} || \le \ \ \ \ $$ $$ \ \ \ \le \frac{B}{(s_1-s_0)^a}\frac{k_0! C^{k_0}}{(1-s_1)^{k_0}} \frac{C'' \ell_1! C^{\ell_1}}{(1-s_1)^{\ell_1}} \dots \ \ \ \ \ \ \  $$ \begin{equation} \ \ \ \ \ \ \ \ \ \ \dots \frac{B}{ (s_{n}-s_{n-1})^a }\frac{k_{n-1}! C^{k_{n-1}}}{(1-s_{n})^{k_{n-1}}} \frac{C''  \ell_n! C^{\ell_n}}{(1-s_{n})^{\ell_n}}  \frac{B}{ (s_{n+1}-s_{n})^a }\frac{k_{n}! C^{k_n}}{(1-s_{n+1})^{k_n}}\frac{m! C^m}{(1-s_{n+1})^m}.  \label{ap24a} \end{equation}

Choose $s_{n+1}$ so that \begin{equation} s_{n+1}-s_n = \frac{a}{m+k_{n}+a} (1-s_n), \label{choice19} \end{equation} then
$$ \frac{1}{(s_{n+1}-s_n)^a(1-s_{n+1})^{k_{n}+m}} = \frac{ (\frac{m+k_{n}+a}{a})^a (\frac{m+k_{n}+a}{m+k_{n}})^{m+k_{n}} }{(1-s_{n})^{k_{n}+m+a}} \le \frac{ (m+k_{n}+a)^a e^a }{a^a (1-s_{n})^{k_{n}+m+a}}. $$
Then
$$ \text{RHS of \eqref{ap24a} \ \ } \le \frac{ (C'' )^n B^{n+1}e^{(n+1)a}}{a^{(n+1)a}(1-s_0)^{(n+1)a}}\left(\frac{ C}{1-s_0}\right)^{k_0+\ell_1+...\ell_n+k_n+m} k_0! \ell_1! k_1!... \ell_n! k_n! m! \times \ \ \ \ $$
$$ \ \ \ \ \ \times (k_n+m+a)^a(k_n+m+\ell_n+k_{n-1}+2a)^a... (k_0+...+m+[n+1]a)^a. $$

But \small
$$ 
 \sum_{\tiny\begin{array}{c} k_0+...+k_n+\ell_1+...+\ell_n+m=N; \\ k_0,...,k_n, \ell_1,..., \ell_n, m\ge 0\end{array}} k_0! \ell_1! k_1!... \ell_n! k_n! m!  (k_n+m+a)^a(k_n+m+\ell_n+k_{n-1}+2a)^a... (k_0+...+m+[n+1]a)^a \ \le \ $$
\begin{equation} \ \ \ \ \ \le (36a^{a-1})^{n+1}(N+[n+1]a)!. \ \ \ \ \label{ap24b} \end{equation}
\normalsize

Indeed, on the LHS of \eqref{ap24b} let us introduce new summation indices:  $m+k_n=:K_n\ge 0$, $\ell_n+k_{n-1} =:K_{n-1}\ge 0$, .... , $\ell_1+k_0=:K_0\ge 0$ and use 
\eqref{fe15in6}; let us also use that $\frac{(p+a)^a}{(p+1)...(p+a)}\le a^{a-1}$ for any $p\ge 0$. Then 
$$ \text{LHS of \eqref{ap24b}} \ \le \ \ \ \ \ \ \ \  \ \ \ \ \ $$
\small
$$ \le (6a^{a-1})^{n+1} \sum_{\tiny \begin{array}{c} K_0+...+K_n=N; \\ K_0,...,K_{n-1},  K_n\ge 0\end{array}} K_0!... K_{n-1}! K_n!   \frac{(K_n+a)!}{K!} \cdot \frac{(K_{n-1}+K_{n}+2a)!}{(K_{n-1}+K_{n}+a)!}... \frac{ (K_0+...+K_n+[n+1]a)!}{(K_0+...+K_n+na)!}   \ = \ $$
$$ \le (6a^{a-1})^{n+1}  \cdot \sum_{\footnotesize \begin{array}{c} K_0+...+K_n=N; \\ K_0,...,K_{n}\ge 0\end{array}} K_0!... K_{n-2}! K_{n-1}!(K_n+a)!  \frac{(K_{n-1}+K_{n}+2a)!}{(K_{n-1}+K_{n}+a)!}... \frac{ (K_0+...+K_n+[n+1]a)!}{(K_0+...+K_n+na)!} \le $$
(put $K'_{n-1}:=K_{n-1}+K_{n}$ and use \eqref{fe15in6} )
\footnotesize
$$ \le (6a^{a-1})^{n+1}  \cdot 6 \sum_{\footnotesize \begin{array}{c} K_0+...+K_{n-2}+K'_{n-1}=N; \\ K_0,...,K_{n-2}, K'_{n-1}\ge 0\end{array}} K_0!... K_{n-2}! (K'_{n-1}+2a)! \frac{  (K'_{n-1}+K_{n-2}+3a)!}{(K'_{n-1}+K_{n-2}+2a)!}...\frac{ (K_0+...+K'_{n-2}+[n+1]a)!}{ (K_0+...+K'_{n-2}+na)!}  \le $$
\normalsize
(etc.)
$$ \le (6a^{a-1})^{n+1} \cdot 6^n (N+[n+1]a)! < (36a^{a-1})^{n+1} (N+[n+1]a)!, $$ 
hence \eqref{ap24b}.

\vskip2pc

Let us go back to interpreting \eqref{ap24c}. Formally writing
$$ \text{``} \ \  \sum_{n\ge 0} h^{na} \calG(\calH\calG)^n v = \sum_{J\ge 0} h^J u_J \ \ \text{"}$$
motivates the definition
$$ u_J := \sum_{na+N=J,\, n,N\ge 0} \sum_{k_0+...+\ell_n+m=N}  \calG_{k_0,s_1,s_0} \calH_{\ell_1,s_1} \calG_{k_1,s_2,s_1} ... \calH_{\ell_n,s_n} \calG_{k_n,s_{n+1},s_n} v_{m,s_{n+1}} ,$$
where, using \eqref{ap24b}
$$ ||u_J|| \le \sum_{na+N=J,\, n,N\ge 0}\frac{(C'' )^n B^{n+1} e^{(n+1)a}}{a^{(n+1)a}}  \frac{C^N}{(1-s_0)^{N+(n+1)a}} (36a^{a-1})^{n+1} (N+(n+1)a)! < $$
(assuming $C''\ge 1$, $C\ge 36a^{a-2}BC''e^a$,  $C\ge 1$ and taking into account that we have $\le J+1<J+a+1$ summands) 
$$  < \frac{(J+a+1)! C^{J+1}}{(1-s_0)^{J+a}} < \frac{(J+a+1)! C^{J++a+1}}{(1-s_0)^{J+a+1}}. $$

\underline{Step 2.} (Proof of \eqref{EvilConclusion})

Note that estimates on $v_{n,s}$ and $u_{n,s}$ for $s$ outside of the range prescribed by $\frac{|h|C}{1-s}<1$ will be used in an essential way.

In the expression \small
\begin{equation} (\calF + h^a \calH) [ \sum_{n\ge 0} (-1)^n h^{an}\sum_{\tiny \begin{array}{c} k_0,...,k_n\ge 0 \\   \ell_1,...,\ell_n,m\ge 0\end{array}} h^{k_0+...+k_n+\ell_1+...+\ell_n+m} \calG_{k_0,s_1,s_0} \calH_{\ell_1,s_1} \calG_{k_1,s_2,s_1} ... \calH_{\ell_n,s_n} \calG_{k_n,s_{n+1},s_n} v_{m,s_{n+1}} ]
\label{ma19su} \end{equation} \normalsize
a well-defined operator is applied to an absolutely convergent sum, once $s_j$'s $j\ge 1$, are chosen separately in each summand as in \eqref{choice19} while $s_0=s$ ; this absolute convergence justifies the manipulations with infinite sums below.

Let $\calI$ denote the set of all indices $n,k_0,..,m$ as in the sum \eqref{ma19su}; to stress dependence of choices of $s_1,...,s_{n+1}$ on $\iota\in \calI$ we will write $s_1(\iota),...,s_{n+1}(\iota)$. The entries of $\iota$ will be denoted $n(\iota), k_0(\iota),..,m(\iota)$; $\sigma(\iota)$ will denote $k_0+...+\ell_n+m$.

Rewrite \eqref{ma19su}:
$$ \text{\eqref{ma19su}} = I + II, $$
where 
$$ I \ = \    \sum_{\footnotesize \begin{array}{c} n\ge 0; \ \ j\ge 0; \ \  k_0,...,k_n\ge 0 \\   \ell_1,...,\ell_n,m\ge 0\end{array}} (-1)^n  h^{an + j + \sigma(\iota)} \calF_{j,s} \calG_{k_0,s_1(\iota),s_0} \calH_{\ell_1,s_1(\iota)} \calG_{k_1,s_2(\iota),s_1(\iota)} \circ ... \ \ \ \ \ \ $$ $$ \ \ \ \ \ \ \ \ \ \ \ \ \ \ \ ... \circ  \calH_{\ell_n,s_n(\iota)} \calG_{k_n,s_{n+1}(\iota),s_n(\iota)} v_{m,s_{n+1}(\iota)};  $$ 
$$ II =    \sum_{\footnotesize \begin{array}{c} n\ge 0; \ \ j\ge 0; \ \  k_0,...,k_n\ge 0 \\   \ell_1,...,\ell_n,m\ge 0\end{array}} (-1)^n h^{a(n+1) + j+\sigma(\iota)} \calH_{j,s} \calG_{k_0,s_1(\iota),s_0} \calH_{\ell_1,s_1(\iota)} \calG_{k_1,s_2(\iota),s_1(\iota)} \circ ... \ \ \ \ \ \ $$ $$ \ \ \ \ \ \ \ \ \ \ \ \ \ \ \ ... \circ   \calH_{\ell_n,s_n(\iota)} \calG_{k_n,s_{n+1}(\iota),s_n(\iota)} v_{m,s_{n+1}(\iota)} .
$$

In $I$ introduce a new index $p=j+k_0$:
 $$ I =  \sum_{\footnotesize \begin{array}{c} n\ge 0; \ \ p\ge 0; \ \  k_1,...,k_n\ge 0 \\   \ell_1,...,\ell_n,m\ge 0\end{array}} (-1)^n  h^{an + p+k_1+...+\ell_n+m} \left(  \sum_{j,k_0\ge 0; \ j+k_0=p} \calF_{j,s} \calG_{k_0,S_1,s_0} \right)\circ \ \ \ \ \ \ $$ $$ \ \ \ \ \ \ \ \ \ \ \ \ \ \ \ \circ  \calH_{\ell_1,S_1} \calG_{k_1,s_2(\iota),S_1} ... \calH_{\ell_n,s_n(\iota)} \calG_{k_n,s_{n+1}(\iota),s_n(\iota)} v_{m,s_{n+1}(\iota)}, $$  
where $S_1$ depends on $n,k_1,...,k_n,\ell_1,...,\ell_n,m,p$ and is defined as $\min \{ s_1(\iota) \} $ over all $\iota$s with prescribed $n,k_1,...,k_n,\ell_1,...,\ell_n,m$ and $k_0\le p$ (so it is a minimum over a set of $p+1$ elements). 
By \eqref{ma20e2} the sum in parentheses is zero unless $p=0$ in which case it is the inclusion $\calY_{S_1}\stackrel{\id}{\hookrightarrow} \calY_{s_0}$. Hence
\begin{equation} I =  \sum_{\footnotesize \begin{array}{c} n\ge 0;  \ \  k_1,...,k_n\ge 0 \\   \ell_1,...,\ell_n,m\ge 0\end{array}} (-1)^n  h^{an +k_1+...+\ell_n+m}  \calH_{\ell_1,s} \calG_{k_1,s_2(\iota),s} ... \calH_{\ell_n,s_n(\iota)} \calG_{k_n,s_{n+1}(\iota),s_n(\iota)} v_{m,s_{n+1}(\iota)} .\label{ma19tar} \end{equation} \normalsize

We see that most terms of $I$ are canceled by terms of $II$, and only the $n=0$ term of  \eqref{ma19tar} remains, which is just $v$ which implies \eqref{EvilConclusion}.
$\Box$

\normalsize 


\section{Proof of the main result} \label{MainPf}

In this section we prove our main result Theorem \ref{MAINTH}.



The main step in the proof of theorem \ref{MAINTH} is, by section \ref{MthOrderFormal}, the solution of equation $F_s(E_0,...E_{M-2},T)=0$, cf. \eqref{ma17f2}, where $F_s: \calX^\tau_s\to \calY^\tau_s$, $s\le 1$, $\tau<\tau_0$, see section \ref{FuncAn}.  With respect to $s$, $F_s$ form a family of maps compatible with inclusions, in the sense of section \ref{GcalcDef} -- for various maps below we will keep in mind and use compatibility with inclusions  without always writing these words. The operator $h$ of section \ref{GcalcDef} will be just the multiplication by $h$ in the spaces $\calX^\tau_s, \calY^\tau_s, \calZ^\tau_s$, with $|h|=\tau$. We will use $h$-linearity properties of various operators without explicitly mentioning it. 

We will take $x^{(0)}\in \calX^\tau_1$ to be the tuple $(E_0,...,E_{M-2},T)$ of {\it polynomials in $h$} that solves $F(E_0,...,E_{M-2},T)$ up to order $h^{b_0+1}$, $b_0=8$, see \eqref{bj},  on the disc of radius $\rho_0$ (see Notation \ref{tau0rho0}); such a tuple exists by Proposition \ref{formalSolu}.

If we write down the expression for $F_s(E_0,...,E_{M-2},T)$ as a power series in $(E_0,...,E_{M-2},T)\in \calX^\tau_s$, the series will converge provided $\| hT'\|_{\sfV^t(s)}<1$ which is definitely the case if \mbox{$\tau\, \| (E_0,...,E_{M-2},T)\|_{\calX^\tau_s} < 1$.}
Therefore, if we impose on $\tau$ the the condition
\begin{equation} \tau \, \| x^{(0)}\|_{\calX^\tau_1} < \frac{1}{2}, \label{DenomFarFr0} \end{equation}
then $F_s(x^{(0)}+w)$ can be represented in the form \eqref{fe12d} with  some finite $\alpha$ independent of $s$.

Since $F_s(E_0,...,E_{M-2},T)$ is an analytic function of $E_0,...,E_{M-2},T$, we can write its 
Fr\'echet derivative by means of usual Calculus formulas. Since $F_s$ explicitly depends only on the $E_0,...,E_{M-2},T,...,T^{(3)}$, we can write $dF_{(E_0,...,E_{M-2},T)}(\calE_0,...,\calE_{M-2},\calT)$ as a finite sum (index $s$ suppressed)
\begin{equation} dF_{(E_0,...,E_{M-2},T)} (\calE_0,...,\calE_{M-2},\calT) = \sum_{\nu=1}^N \calF_{\nu}(E_0,...,E_{M-2},T) \cdot L_\nu(E_0,...,E_{M-2},\calT),  \label{FinLinCom}\end{equation} 
where $(\calE_0,...,\calE_{M-2},\calT)$ belong to the tangent space of $\calX^\tau_s$ at $(E_0,...,E_{M-2},T)$, and $\calF_\nu:\calX^\tau_s\to \calY^\tau_s$ are analytic on an open subset of $\calX^\tau_s$ compatible with inclusions $\calX^\tau_s\to \calX^\tau_{s'}$ and $\calY^\tau_s\to \calY^\tau_{s'}$, and $L_\nu:\calX^\tau_s\to \calY^\tau_s$ are constant (i.e. $(E_0,...,E_{M-2},T)$-independent) linear maps $\calX^\tau_s\to\calY^\tau_s$ of norm $\le 1$, also compatible with inclusions. The dot $(\cdot)$ on the RHS of \eqref{FinLinCom} is the pointwise multiplication of elements of $\calY^\tau_s$ which  has norm $1$ in the sense that $\| y_1 \cdot y_2\|_{\calY^\tau_s}\le \| y_1\|_{\calY^\tau_s} \, \| y_2\|_{\calY^\tau_s}$.

For concreteness, variating the LHS  \eqref{ma17f2}, we pick $N=M+5$ and make the following choices where specifics of big formulas will be unimportant later on: 
\footnotesize
\begin{equation} \begin{array}{lllcl} \calF_1 = &&  1, && L_1 = \sum_{\mu=0}^{m-2} z^{\mu}\calE_{\mu} -\frac{z^M}{2} \calT' - \frac{M z^{M-1}}{4} \calT; \\
\calF_{2+j}=& h\hat\calF_{2+j} =&  \sum_{\mu=1}^{j} \binom{j}{\mu} h^\mu T^\mu z^{j-\mu}  +  h T' \left( 2 + hT' \right)   (z+hT)^j   && L_{2+j} = \calE_j, \ \ \ j=0,...,M-2; \\
\calF_{M+1}= & h\hat\calF_{M+1} & \text{(see below)}
 && L_{M+1} = \calT \\
\calF_{M+2} = & h\hat \calF_{M+2}  & \text{(see below)} 
 && L_{M+2}=\calT' 
\\
\calF_{M+3} = & h^2 \hat{\hat \calF}_{M+3} = &
\frac{h^3}{2}  \frac{T'''}{(1+hT')^2 }-  \frac{3h^4}{2}  \frac{ {T''}^2}{( 1+hT' )^{3}}   && L_{M+3}= \calT' 
\\
\calF_{M+4}= & h^2 \hat{\hat \calF}_{M+4} = & \frac{3h^3}{2}  \frac{ {T''}}{( 1+hT' )^{2}}
 && L_{M+4}=\calT''
\\
\calF_{M+5} = & h^2 \hat{\hat \calF}_{M+5} = & -  \frac{h^2}{2}  \frac{1}{1+hT'} && L_{M+5} = \calT''',
\end{array}
\label{ma30a} \end{equation}
\normalsize
where \footnotesize
$$ \hat\calF_{M+1} =   (1+hT')^2 \sum_{j=0}^{M-2} j E_j  (z+hT)^{j-1}  - \frac{\sum_{\mu=2}^M \binom{M}{\mu} h^{\mu-2} z^{M-\mu}  \mu  T^{\mu-1}}{4} - \frac{ M T'\left( 2 + hT' \right)  (z+hT)^{M-1}}{4}, $$
$$ \hat \calF_{M+2} =  2  \left( 1 + hT' \right) \left( \sum_{j=0}^{M-2} E_j  (z+hT)^j - \frac{1}{4} \sum_{\mu=1}^M z^{M-\mu} h^{\mu-1} T^{\mu}\right). $$
\normalsize

Notice that norms of $L_\nu$ are all $\le 1$ thanks to the assumption on $\rho_0$ that we made in Notation \ref{tau0rho0}. \label{ma31lab5pm}


Motivated by our considerations in section \ref{dF0}, we have grouped the terms $\calF_\nu L_\nu$ in such a way that terms $\nu=1,..., M+2 $ only depend on $E_0,...,E_{M-2},T,T'$ and $\calE_0,...,E_{M-2},\calT,\calT'$; terms containing higher derivatives are put into summands for $\nu=M+3,M+4,M+5$, we notice that all summands the latter group also contains a factor of $h^2$. We will thus treat $\calF_{M+3}L_{M+3}+\calF_{M+4}L_{M+4}+\calF_{M+5}L_{M+5}$ as a perturbation of $\calF_1 L_1+...+\calF_{M+2}L_{M+2}$.
In turn, as $\calF_2,...,\calF_{M+2}$ contain a prefactor of $h$, we treat $\calF_2L_2+...+\calF_{M+2} L_{M+2}$ as a perturbation of $\calF_1 L_1$.

Modify $\alpha$ in such a way that representations of the form \eqref{fe12d} also hold for $\hat\calF_1(x^{(0)}+w),...,\hat\calF_{M+2}(x^{(0)}+w),\hat{\hat\calF}_{M+3}(x^{(0)}+w),...,\hat{\hat\calF}_{M+5}(x^{(0)}+w) $. 

%

From now on we fix the following notation:\\
$A$ is the constant such that $L_1$ understood as a map $\calZ^\tau_s \to \calY^\tau_s$ has an inverse of norm $\le A$;\\
$B$ is the constant such that the restriction map $\calZ^\tau_s\to \calX^\tau_{s'}$, $s'<s$, is of norm $\le \frac{B}{(s-s')^2}$; \\
we fix the following functions:  
$$ \beta(t) = te^{4t}; \ \ \ \ \beta_2(t) =\sum_{k\ge 1} \frac{t^{k+1} 4^k}{k!}. $$ 

We will need an integer sequence $b_j$, $j\ge 0$, satisfying the properties:
\begin{equation} \text{a)} \ b_{j+1} =2b_j-7; \ \ \ \text{b)} \ b_j\ge 3+j ;
\label{bj}
\end{equation}
clearly the condition a) together with $b_0=8$ generates such a sequence.

We assume that $\tau$ satisfies \eqref{DenomFarFr0}, \eqref{y0cond}, \eqref{hm12},  \eqref{ma30hn},  \eqref{Hless12}, and \eqref{ma16ni}; we assume that $\alpha$ is a finite number chosen as above, and we assume that $C$ satisfies \eqref{m16e} and \eqref{y0cond}. 

\underline{Step 1:} Constructing iterations of the Newton method as formal factorially divergent expansions.

We will construct Gevrey expansions which will constitute iterations of the Newton method.  For our intuition, we suggest the following correspondence (which we will not make precise) of the notions of section \ref{NMsec} to the objects introduced below:   
$$ x^{(j)}_{\text{sec.\ref{NMsec}}} \leftrightarrow x^{(0)}+g^{(j)}, $$
$$ w^{(j)}_{\text{sec.\ref{NMsec}}} \leftrightarrow \sum_{n\ge 1} h^n \tilde w^{(j)}_n \ \ \text{and} \ \ \ h^{b_j-3} \sum_{n\ge 1} h^n w^{(j)}_n, $$
$$ y^{(j)}_{\text{sec.\ref{NMsec}}} \leftrightarrow h^{b_j} \sum_{n\ge 0} h^n y_n^{(j)}. $$

Speaking rigorously, we will now construct the following objects for all $j\ge 0$:
\begin{equation} \text{for} \ n\ge 0 \ \ y^{(j)}_n = (y^{(j)}_{n,s})\in \mathring \calY^\tau_1; \ \ \| y^{(j)}_{n,s}\| \le \frac{n! C^n}{(1-s)^n}; \label{m16ind} \end{equation}
\begin{equation} \text{for} \ n\ge 1 \ \ \tilde w^{(j)}_n = (\tilde w^{(j)}_{n,s})\in \mathring \calX^\tau_1; \ \ \ \ 
\| \tilde w^{(j)}_{n,s} \| \le \frac{1}{2^j} \frac{C^n n!}{(1-s)^n}. \label{m16W} \end{equation}

We will use the notation:
\begin{equation} \text{for} \ n\ge 1 \ \  g^{(0)}_{n,s}=0; \ \ \ g^{(j)}_{n,s} = \tilde w^{(0)}_{n,s}+...+\tilde w^{(j-1)}_{n,s}, \ \ \ j\ge 1 ; \label{m16G} \end{equation}
clearly $g^{(j)}_{n} \in \mathring X^\tau_1$ and, as soon as \eqref{m16W} is true, $g^{(j)}_{n,s}$ satisfy  \eqref{ma6a} with $C_0=2$, i.e.  
$$ \| g_{n,s}^{(j)} \|_{\calX^\tau_s} \le 2 \cdot \frac{C^n n!}{(1-s)^n} .$$ 

By our choice of $x^{(0)}$, we can write $F(x^{(0)})$ as a convergent power series $h^{b_0+1} \sum_{n\ge 0} \tilde y_n(z) h^n$, where $\tilde y_n(z)$ is an $h$-independent function and $\| \tilde y_n(z) \|_{\calY^\tau_1} \le c' \cdot c^n$. Put $y^{(0)}_{n,s} := h \tilde y_n(z)$, then assuming 
\begin{equation} C\ge c \ \ \ \text{and} \ \ \ \tau \cdot c' \le 1 \label{y0cond} \end{equation}
we obtain \eqref{m16ind} for $j=0$.

Suppose $y^{(j')}_{n,s}$ are defined for $j'\le j$ and $w^{(j')}_{n,s}$ are defined for $j'<j-1$ (vacuously true if $j=0$). Let us construct $w^{(j)}$ and $y^{(j+1)}$.

Using lemma \ref{InsertL} with $x^{(0)}$ for $x_0$, with $g^{(j)}_{n,s}$ for $g_{n,s}$,  
and with $\hat \calF_\nu$ and $\hat {\hat \calF}_\nu$ for $f$, delivers elements 
\begin{equation}  G_{\nu,n} \in \mathring \calY^\tau_1,  \ \ \ ||G_{\nu,n,s}||_{\calY^\tau_s} \le \beta(C_0\alpha) \, \frac{C^n n!}{(1-s)^n}, \ \ \ \ n\ge 1, \ \  \nu=2,...,M+2; \label{m16a} \end{equation}
\begin{equation}   G'_{\nu,n} \in \mathring \calY^\tau_1, 
 \ \ \| G'_{\nu,n,s} \|_{\calY^\tau_s} \le \max\{ \| \hat{\hat\calF}_\nu(x^{(0)})\|_{\calY^\tau_1} , \beta(C_0\alpha)\} \frac{C^n n!}{(1-s)^n}, \ \ \  n\ge 0, \ \ \nu=M+3,M+4,M+5, \label{m16b} \end{equation}
such that, intuitively speaking, 
$ \hat \calF_\nu(x^{(0)}+\sum_n h^n g^{(j)}_n)$ corresponds to $\hat \calF(x^{(0)})+ \sum_{n\ge 1} h^n G_{\nu,n}$, and $\hat{\hat\calF}_\nu(x^{(0)}+\sum_n h^n g^{(j)}_n)$ corresponds to $\sum_{n\ge 0} h^n \calG'_{\nu,n}$.

Let us find an infinite series of operators $\calL_{n,s}:\calY^\tau_s \to \calZ^\tau_s$ so that the infinite sum $\sum_{n\ge 0} h^n \calL_{n,s}$ will play the role of the inverse of $\sum_{\nu=1}^{M+2} \calF_\nu(x) \cdot L_\nu$ in the precise sense specified below. An estimate  of $\| \calL_{n,s}\|/ \left( \frac{C^n n!}{(1-s)^n} \right)$  by $2A$ will be important in \eqref{Ble} which, in turn, lets us preserve the same constant $C$ from one induction step to the other. 

With this goal in mind, we apply the lemma \ref{BenignInv} with operators
\begin{equation} (\calF)_{\rm Lemma} = \sum_{\nu=1}^4 \calF_\nu(x^{(0)})\cdot L_\nu, \ \ \ (\calH_{n,s})_{\rm Lemma} :=  \sum_{\nu=2}^{M+2} h G_{\nu,n,s}\cdot L_\nu , \ \ \ \  n\ge 1. \label{ma15e1} \end{equation}
and constants
$$ (C_1)_{\rm Lemma}=(M+1)\tau  \beta(C_0\alpha); \ \ \ A_{\rm Lemma} = 2A_{\rm here}. $$
In order to assure that, as the assumptions of lemma \ref{BenignInv} require, $\calF$ has inverse of norm $\le 2A$, we remember that  $L_1$ has an inverse of norm $\le A$; it is for that reason that we assumed that $\tau$ should satisfy
\begin{equation} \tau \cdot \sum_{\nu=2}^{M+2} \| \hat \calF_\nu(x^{(0)})\|_{\calY^\tau_1} \le \frac{1}{A}. \label{hm12} \end{equation}

We have assumed above that  $\tau$ is so small that 
\begin{equation}   ((C_1)_{\rm Lemma} A)^2 e^{4 (C_1)_{\rm Lemma} A}\ \le \ 2A, \label{ma12hn} \end{equation}
or, more explicitly,
\begin{equation} ((M+1)\tau A \beta(2\alpha) )^2 e^{4 (M+1)\tau  A \beta(2\alpha) }\ \le \ 2A; \label{ma30hn}  \end{equation}
therefore, the lemma \ref{BenignInv} yields operators
$$ \calL_{n,s} : \calY^\tau_s\to \calZ^\tau_s \ \ \ \text{s.th.} \ \ \| \calL_{n,s}\|\le 2A \frac{C^n n!}{(1-s)^n} $$
satisfying \eqref{ma20e1}.

For any $s'<s$ we can compose the operator $\calL_{n,s}$ with the restriction map $\calZ^\tau_s\to \calX^\tau_{s'}$  and get $\calL'_{n,s,s'}: \calY^\tau_s \to \calX^\tau_{s'}$ of norm $\le \frac{B}{(s-s')^2} \cdot 2A \frac{C^n n!}{(1-s)^n}$.



We are thus in the situation of Lemma \ref{EvilInv} with $a=2$, operators (cf. \eqref{m16a}, \eqref{m16b})
\begin{equation} (\calF_0)_{\rm Lemma} = \sum_{\nu=1}^{M+2} \calF_\nu(x^{(0)})\cdot L_\nu, \ \ \  (\calF_n)_{\rm Lemma} = \sum_{\nu=2}^{M+2} G_{\nu,n}\cdot L_\nu, \ \ \ \text{if} \ n\ge 1;  \label{m27b} \end{equation} \begin{equation}  (\calH)_{\rm Lemma} = \text{``} \sum_{\nu=M+3}^{M+5} \hat{\hat\calF}_\nu(x) \cdot L_\nu \text{"}, \ \ \ \text{i.e.} \ \ (\calH_{n})_{\rm Lemma}= \sum_{\nu=M+3}^{M+5} G'_{\nu,n} \cdot L_\nu, \label{m16c} \end{equation} 
vectors (see \eqref{m16ind}\ ) 
\begin{equation} (v_n)_{\rm Lemma} = - y^{(j)}_n, \ \ \ \ n\ge 0; \label{m27d} \end{equation}
and constants
$$  (C')_{\rm Lemma} = \max\{ \sum_{\nu=1}^{M+2} |\calF_\nu(x^{(0)})|_{\calX^\tau_1} , \, C_1 \}; $$ \begin{equation} (C'')_{\rm Lemma}=C'':=\sum_{\nu=M+3}^{M+5} \max\{ \| \hat{\hat\calF}_\nu(x^{(0)})\|_{\calY^\tau_1}, \beta(C_0\alpha)\};    
  \label{Ble} \end{equation}
$$ B_{\rm Lemma} = 2 B_{\rm here} A, $$
and we must assume 
\begin{equation} C\ge \max\{ 1, 36 \cdot 2BA  \cdot C'' \cdot e^2 \}. \label{m16e} \end{equation} 
Then Lemma \ref{EvilInv} delivers vectors $u_n$ s.th. 
$$  || u_{n,s}||_{\calX^\tau_s} \le \frac{C^{n+3} (n+3)!}{(1-s)^{n+3}} $$
Put $w^{(j)}_n=u_{n-3}$ for $n\ge 3$, $w^{(j)}_1=w^{(j)}_2=0$, then $w=\sum_{n\ge 1} h^n w_n$ plays the role of the solution of 
$$ (dF)_x(w) \ = \ h^3 \sum_{n\ge 1} h^n y^{(j)}_n. $$
Define $\tilde w^{(j)} = h^{b_j-3} w$. Because of \eqref{bj}, b) and  because have assumed 
\begin{equation}   \tau<\frac{1}{2}, \label{Hless12} \end{equation}
we see that the estimate of \eqref{m16W} is satisfied.

%
%

\vskip2.5pc

We will now  obtain vectors $y^{(j+1)}_n$ such that, intuitively speaking, $h^{b_{j+1}}\sum_{n\ge 1} h^n y^{(j+1)}_n$ plays the role of $f(x^{(0)}+g^{(j)}+w^{(j)})-f(x^{(0)}+g^{(j)})-[f'(x^{(0)}+g^{(j)})(w^{(j)})]$ (the cancellation of the second and third summands will be discussed later). apply lemma \ref{DiffWellDef},b) with the following inputs: 
$$ f_{\text{Lemma}} = F_{\rm here};  \ \ \ (x_0)_{\rm Lemma} = x^{(0)}, \ \ \ \ g_{\text{Lemma}} = g^{(j)}, $$
$$ b_{\rm Lemma} = b_j-3 , \ \ \ (v_n)_{\rm Lemma} = w_n^{(j)}, $$ 
$$   \alpha_{\text{Lemma}} = \alpha_{\rm here} , \ \ \ (C_0)_{\text{Lemma}} = 2 , \ \ \ A_{\text{Lemma}} = 1. $$

Then Lemma \ref{DiffWellDef},b) yields
$$ u_{n,s}\in \calY^\tau_s, \ \ n\ge 1$$
satisfying 
$$ \| u_{n,s}\|_{\calY^\tau_s} \le \beta_2(\alpha_{\rm Lemma} ( (C_0)_{\rm Lemma} + A_{\rm Lemma}) \frac{n! C^n}{(1-s)^n} =  \beta_2(3\alpha) \frac{n! C^n}{(1-s)^n}$$
such that $h^{2b_j-6} \sum_{j\ge 1} h^j u_j$ plays the role of $y^{(j+1)}$. 

In \eqref{bj},a) we have defined $b_{j+1} = 2b_j - 7$; now take  $y^{(j+1)}_n = h u_n$ and assume that 
\begin{equation} \tau \cdot \beta_2(3\alpha) \le 1. \label{ma16ni}\end{equation} 
Then 
$y^{(j+1)}_n$ satisfy the inductive assumption \eqref{m16ind}. 

The inductive construction of step 1 is thus complete.

\vskip2.5pc

\underline{Step 2.} (Passing from Gevrey expansions to actual vectors)

If condition 
\begin{equation} \tau C \ < \ 1-s\label{GevConv} \end{equation} is satisfied, then we can define the following vectors by means of absolutely convergent sums:
\begin{equation} \bfw^{(j)} = \sum_{n\ge 1} h^n \tilde w^{(j)}_{n,s}  \stackrel{\text{abs.conv.}}{=} h^{b_j-3} \sum_{n\ge 1} h^n w^{(j)}_{n,s} \ \in \ \calX^\tau_s ; \label{bfwDef} \end{equation}
\begin{equation} \bfg^{(j)} = \bfw^{(0)} +...+\bfw^{(j-1)} \stackrel{\text{abs.conv.}}{=} \sum_{n\ge 1} h^n g^{(j)}_{n,s} \ \in \ \calX^\tau_s ; \label{bfgDef} \end{equation}
\begin{equation} \bfy^{(j)} = h^{b_j} \sum_{n\ge 0} h^n y_n^{(j)} \ \in \ \calY^\tau_s; \label{bfyDef} \end{equation} 
the equalities marked above as (abs.conv.) are justifiable by operations on absolutely convergent sums.

\underline{Claim.} If \eqref{cn2m27} and \eqref{co27m} are satisfied, then $ F(x^{(0)}+\bfg^{(j)}) = \bfy^{(j)}$. 

The case $j=0$ is obvious by definitions.

Suppose the Claim is true for $j$, let us deduce it for $j+1$.

Indeed, applying Lemma \ref{DiffWellDef} with the same ingredients as in Step 1 but under condition 
\begin{equation} \sum_{n=1}^\infty \frac{\tau^n C^n}{(1-s)^n} 
< \frac{1}{3\alpha}
\label{cn2m27} \end{equation} (which we have assumed above) corresponding to \eqref{ma8c}, we have
$$ F(x^{(0)}+\bfg^{(j)}+w^{(j)}) = \bfy^{(j+1)} + F(x^{(0)}+\bfg^{(j)}) + [F'(x^{(0)}+\bfg^{(j)})] (\bfw^{(j)}) $$
so our claim, in view of induction hypothesis, is reduced to showing that 
\begin{equation} \bfy^{(j)} \ = \ -[F'(x^{(0)}+\bfg^{(j)})] (\bfw^{(j)}) \label{ma27a} \end{equation}

By Lemma \ref{EvilInv}, $\bfw^{(j)}$ solves the equation
$$ (\bfF + h^2 \bfH) \bfw^{(j)} = -\bfy^{(j)} $$
where (cf. \eqref{m27b}, \eqref{m16c} )
$$\bfF = \sum_{\nu=1}^4 \calF_\nu(x^{(0)})\cdot L_\nu + \sum_{n\ge 1}  h^n \sum_{\nu=2}^4 G_{\nu,n,s}\cdot L_\nu;  $$ 
$$\bfH = \sum_{n\ge 0} h^n \sum_{\nu=5}^7 G'_{\nu,n,s}\cdot L_\nu.$$ 

But the definition of $G_\nu$ and $G'_\nu$ by means of lemma \ref{InsertL} shows that under  condition 
\begin{equation} \sum_{n\ge 1} \frac{\tau^n C^n}{(1-s)^n} < \frac{1}{2\alpha} \label{co27m} \end{equation} (which we have assumed)
corresponding to \eqref{con11f}, we have
$$ \bfF = \sum_{\nu=1}^4 F(x^{(0)}+\bfg) ; \ \ \  \sum_{\nu=5}^7 \hat{\hat\calF}_\nu(x^{(0)}+\bfg) \cdot L_\nu  = \bfH $$
and so by \eqref{FinLinCom} we conclude the proof of the claim.

\underline{Step 3.} (Proving the convergence of the Newton's method.)

%
By the above construction, 
$$ \| \bfy^{(j)} \|_{\calY^\tau_s} \ \le \ \tau^{b_j} \sum_{n\ge 0} \frac{\tau^n C^n}{(1-s)^n} < \frac{1}{2^j}\left(1+\frac{1}{3\alpha}\right),  $$
$$ \| \bfw^{(j)} \|_{\calX^\tau_s} \ \le \ \frac{1}{2^j} \sum_{n\ge 1} \frac{\tau^n C^n}{(1-s)^n} \le \frac{1}{3\alpha 2^j}. $$
Thus $x^{(0)}+\bfg^{(j)}$ has a limit $x^{(0)}+\bfg^{\infty}$ in $\calX^\tau_s$ and $F(x^{(0)}+\bfg^{\infty})=0$. 

Let us review the choices of various constants that we have made.  In the beginning of this section we have chosen and fixed $\alpha$ large enough depending on the initial data (namely, on $\tilde Q_1(z,h)$). In Step 1 we have assumed that $\tau$ is small enough to satisfy  \eqref{DenomFarFr0}, \eqref{y0cond}, \eqref{hm12},  \eqref{ma30hn},  \eqref{Hless12}, and \eqref{ma16ni} (which depend only the initial data); then we have chosen and fixed $C$ large enough to satisfy \eqref{m16e} and \eqref{y0cond}; a choice of $C$ valid for one value of $\tau$ is also valid for smaller values of $\tau$. In Step 2, in \eqref{GevConv} we choose and fix $s$; a choice of $s$ valid for one value of $\tau$ will also be valid for smaller values of $\tau$. Then we shrink $\tau$ if necessary to satisfy \eqref{cn2m27} and \eqref{co27m}.
We conclude that there exist $s_*$ and $\tau_*$ such that for any $s<s_*$ and $\tau<\tau_*$, the equation $F(\bfx) =0$ has a solution $\bfx\in \calX^\tau_s$.

\vskip2.5pc

Finally, in view of section \ref{MthOrderFormal} we have shown that $E_0,...,E_{M-2},T$, and hence also $y(x,h)$ have the Gevrey growth condition. Using \eqref{factorpsi}, one shows by a standard argument that also $\psi(x,h)$ satisfies a Gevrey growth condition. This concludes the proof of theorem \ref{MAINTH}. $\Box$

{\bf \large Acknowledgements.} The author would like to thank Yoshitsugu Takei for his help in reading ~\cite{AKT}, Pierre Schapira for a very clear introduction to scales of Banach spaces in ~\cite[Ch.{\S}I.3]{Sch}, as well as Stavros Garoufalidis, Gleb Novichkov, and Ilya Roublev who found time to answer author's questions.

\end{document}